\documentclass{article}
\usepackage{amsmath,latexsym,amsfonts}

\numberwithin{equation}{section}

\newtheorem{theorem}[equation]{Theorem}
\newtheorem{lemma}[equation]{Lemma}
\newtheorem{corollary}[equation]{Corollary}

\newtheorem{example}[equation]{Example}
\newtheorem{remark}[equation]{Remark}

\newcommand{\Pic}{\operatorname{Pic}}
\newcommand{\Div}{\operatorname{div}}
\newcommand{\mult}{\operatorname{mult}}
\newcommand{\codim}{\operatorname{codim}}
\newcommand{\sing}{\operatorname{Sing}}

\newcommand{\veps}{\varepsilon}
\newcommand{\PP}{{\mathbb P}}
\newcommand{\Q}{{\PP^2}}
\newcommand{\CO}{{\mathcal O}}
\newcommand{\CI}{{\mathcal I}}
\renewcommand{\L}{{\mathcal L}}

\newcommand{\mapright}[1]{\stackrel{#1}\longrightarrow}

\newcommand{\mapdown}[1]{\Big\downarrow {\scriptstyle #1}}
\newcommand{\prox}{\to}

\newcommand{\qed}{\hfill{$\Box$}}

\newenvironment{pf}
{\noindent\textbf{Proof.}}
{\medskip}

\newcommand{\fref}[1]{(\ref{#1})}

\newcommand{\ndef}{{\mathrm{def}}}
\newcommand{\Def}{{\mathrm{Def}}}

\newcounter{citaz}
\newcommand{\newcit}[1]{\refstepcounter{citaz}\label{#1}}

\begin{document}
\title{Explicit Resolutions of Double Point\\ Singularities of Surfaces}
\author{{\sc Alberto Calabri}\thanks{\ Dipartimento di Matematica,
Universit\`{a} di Roma ``Tor Vergata'',
via della Ricerca Scien\-ti\-fi\-ca, 00133 Roma.
E-mail: {\tt calabri@mat.uniroma2.it}.
Partially supported by E.C.\ project EAGER,
contract n.\ HPRN-CT-2000-00099.}
\and {\sc Rita Ferraro}\thanks{\ Dipartimento di Matematica,
Universit\`{a} di Roma Tre,
Largo S.~Leonardo Murialdo 1, 00146 Roma.
E-mail: {\tt ferraro@mat.uniroma3.it}.}
}

\date{}
\maketitle

\begin{abstract}
Locally analytically, any isolated double point occurs as a double
cover of a smooth surface. It can be desingularized explicitly
via the canonical resolution, as it is very well-known.
In this paper we explicitly compute the fundamental cycle
of both the canonical and minimal resolution of a double point singularity
and we classify those for which the fundamental cycle
differs from the fiber cycle.
Moreover we compute the conditions that
a double point singularity imposes to pluricanonical systems.
\end{abstract}



\section{Introduction}


In this article, we give a detailed analysis of isolated surface
singularities of multiplicity two, i.e., \emph{double points}. Our
goal is to be as explicit as possible.

A neighbourhood of an isolated double point $p$ on a complex
(normal) surface is analytically isomorphic to a double cover
$\pi_0:X_0\to Y_0$ branched along a reduced curve $B_0$ with an
isolated singularity at a point $q_1$, where $Y_0$ is a smooth
surface. If $x,y$ are local coordinates of $Y_0$ near $q_1$, then
we may assume that $X_0$ is defined locally by an equation
$z^2=f(x,y)$, where $B_0$ is $f=0$ and $f$ is a square-free
polynomial in $x,y$.

It is very well-known that one may desingularize
$p=\pi_0^{-1}(q_1)\in X$ following the \emph{canonical}
resolution, which consists in desingularizing the branch curve and
normalizing.
This fact was implicitly used by the so-called Italian school
(as one finds out for example by reading Castelnuovo and Enriques' papers
on double covers of the projective plane),
but it has been first written in 1946 by Franchetta \cite{fr1},
who computed the properties of the exceptional curves
of the canonical resolution and the \emph{fiber cycle},
that is the maximal effective divisor contained in the scheme-theoretic
fiber over $p$.

Later, in 1978, Dixon \cite{dixon} proved again Franchetta's results
in a more modern language and he found sufficient conditions
for the \emph{fundamental cycle} of a double point singularity
to be equal to the fiber cycle.

At the same time, Horikawa \cite{ho} and Laufer \cite{laufer3} explained
the canonical resolution process too
(see also \cite{bpv}).
In particular Laufer proved that
the minimal resolution of a double point singularity is
obtained from the canonical one by contracting simultaneously
finitely many
disjoint $(-1)$-curves and he described the relation between
the topological types of the canonical resolution and of the minimal one,
by using essentially the properties of the fundamental cycle.

We recall that the fundamental cycle of a resolution may well be
computed inductively from the knowledge of the intersection matrix
of the exceptional curves. However no explicit formula was known
in the general case. In Theorem \ref{numcycle}, we will give and
prove such a formula, that turns out to be very simple, for the
canonical resolution and then in \fref{barF} for the minimal one.

Moreover, this formula allows us also to classify those double
points for which the fiber cycle strictly contains the fundamental cycle.

Finally in the last section we compute what are the conditions that
a double point singularity imposes to canonical and pluricanonical
systems of a surface.
For this purpose, it is convenient to consider projective surfaces:
we will assume that $X_0$ is a double plane, i.e., $Y_0=\Q$.

Although the general techniques used in this paper are
very well-understood,
we found no reference about our results.
Thus we hope that this paper may serve as a natural complement
to Dixon's and Laufer's papers and as an adequate reference
for algebraic geometers that may use these results.

Moreover our approach is slightly different from the previous ones
because we use, as discrete invariant of a plane curve singularity,
the \emph{Enriques} digraph, namely the directed graph involving
the proximity relations among the infinitely near points to $q_1$,
and we believe that this approach may be of independent interest as well
(cf.\ Enriques diagrams in \cite{kp}, \cite{casas} and \cite{roe}).

Thus we show also some examples of double point singularities
which may help to understand the features of
the combinatorial machinery we introduced.
Please do not hesitate to contact the authors if you want to see
the implementation of this approach on a computer.

The interested reader may also consult \cite{BMcEN}
for the analysis of an embedded resolution of
$\{z^2=f(x,y)\}\subset{\mathbb C}^3$.

\medskip\noindent{\bf Acknowledgements.}
We warmly thank prof.\ Ciro Ciliberto for addressing us to this subject
and prof.\ Rick Miranda for several useful discussions and for joining us
in preparing a preliminary version of this paper.
We are grateful to prof.\ J.~Lipman for some relevant bibliographical reference
and to the referee for many suggestions that really improved the exposition
of this paper.

\section{Notation}

To help the reader with the notation, which will soon
become very heavy, we list here the main symbols used in this paper,
together with the meaning and the reference formula or page
where they are defined.


\noindent
\begin{tabular}{l@{\,\,}l@{}r}
\small\it Symbol & \small\it Meaning & \small\it Ref. \\ \hline \\[-3mm]
$X_0$ & normal complex surface
   \dotfill & p.~\pageref{X0} \\
$p$   & isolated double point singularity of $X_0$
   \dotfill & p.~\pageref{p} \\
$Y_0$ & smooth surface
   \dotfill & p.~\pageref{Y0} \\
$\pi_0:X_0\to Y_0$ & double cover
   \dotfill & p.~\pageref{pi0} \\
$B_0$ & branch curve of $\pi_0$
   \dotfill & p.~\pageref{B0}  \\
$q_1$ & $=\pi_0(p)$, isolated singular point of $B_0$
   \dotfill & p.~\pageref{q1} \\
$\sigma_i:Y_i\to Y_{i-1}$ & blowing-up at a point $q_i\in Y_{i-1}$
   \dotfill & p.~\pageref{sigma_i} \\
$q_i$ & center of the blowing-up $\sigma_i$
   \dotfill & p.~\pageref{q_i} \\
$(\ \cdot\ )_i,\ (\ \cdot\ )$ & intersection pairing in $Y_i$ (resp.\ in $Y_n$)
   \dotfill & p.~\pageref{()i} \\
$\sigma_{ij}:Y_j\to Y_i$ & $=\sigma_j\circ\sigma_{j-1}\circ\cdots\circ\sigma_{i+1}$
   \dotfill & \fref{sigma_ij} \\
$\sigma$ & $=\sigma_{0n}$, sequence of blowing-ups
   \dotfill & p.~\pageref{sigma_ij} \\
$E_i$ & (proper transform of) the exceptional curve $\sigma^{-1}_{i}(q_i)$
   \dotfill & p.~\pageref{defEi} \\
$E^*_i$ & total transform of $E_i$ in $Y=Y_n$ via $\sigma=\sigma_{0n}$
   \dotfill & p.~\pageref{defE*i} \\
$N=(n_{ij})$ & $n\times n$ matrix, $E_i=\sum_{j=1}^n n_{ij}E_j^*$
   \dotfill & \fref{Ei}\\
$M=(m_{ij})$ & $n\times n$ matrix, $E^*_j=\sum_{k=1}^n m_{jk}E_k$
   \dotfill & \fref{Ei}\\
$q_j\to q_i$ & $q_j$ is proximate to $q_i$
   \dotfill & p.~\pageref{qjproxqi} \\
$q_j>^s q_i$ & $q_j$ is infinitely near of order $s$ to $q_i$
   \dotfill & p.~\pageref{qj>qi} \\
$Q=(q_{ij})$ & $n\times n$ matrix, $q_{ij}=1$ if and only if $q_j\to q_i$
   \dotfill & p.~\pageref{defQ}\\
$S=(s_{ij})$ & $n\times n$ matrix of intersection numbers $s_{ij}=(E_i \cdot E_j)$
   \dotfill & p.~\pageref{defS} \\
$\tilde B_i$ & proper transform of $B_0$ in $Y_i$
   \dotfill & p.~\pageref{tildeBi} \\
$\tilde\alpha_i$ & multiplicity of $\tilde B_{i-1}$ at $q_i$
   \dotfill & \fref{tildealfaj} \\
$\tilde\beta_i$ & $=\sum_{j=1}^{n}\tilde\alpha_i m_{ji}$
   \dotfill & \fref{tildebetaj} \\
$\tilde\Gamma_i$ & $=\tilde B_{n|E_i}$, divisor on $E_i$
   \dotfill & p.~\pageref{tildegammaj} \\
$\tilde\gamma_i$ & $=\deg(\tilde\Gamma_i)={\tilde B_n}\cdot{E_i}$
   \dotfill & p.~\pageref{tildeGammaj} \\
$\pi_i:X_i\to Y_i$ & normal double cover induced by $\pi_0$ and $\sigma_{0i}$
   \dotfill & p.~\pageref{pi_i} \\
$B_i$ & branch curve of $\pi_i$
   \dotfill & p.~\pageref{pi_i} \\
$\mu_i$ & multiplicity of $B_{i-1}$ at $q_i$
   \dotfill & p.~\pageref{pi_i} \\
$\veps_i$ & $=\mu_i \bmod 2\in\{0,1\}$, branchedness of $E_i$
   \dotfill & \fref{p=alfa'} \\
$\alpha_i$ & $=\mu_i-\veps_i$, even integer number
   \dotfill & \fref{Balfabeta} \\
$\beta_i$ & $=\sum_{j=1}^{n}\alpha_i m_{ji}$, even integer number
   \dotfill & \fref{Balfabeta} \\
$\tau:X\to X_0$ & canonical resolution of $p\in X_0$
   \dotfill & p.~\pageref{tau} \\
$\bar\tau:\bar X\to X_0$ & minimal resolution of $p\in X_0$
   \dotfill & p.~\pageref{minresthm}  \\
$Y$  & $=Y_n$, smooth surface
   \dotfill & p.~\pageref{Y} \\
$\pi:X\to Y$ & smooth double cover induced by the canonical resol.\ $.\,.$
   \dotfill & p.~\pageref{pi:} \\
$B$ & $=B_n$, smooth branch curve of $\pi$
   \dotfill & p.~\pageref{defB} \\
$\Gamma_i$ & $=B_{|E_i}$, divisor on $E_i$
   \dotfill & \fref{Wi} \\
$\gamma_i$ & $=\deg(\Gamma_i)={B}\cdot{E_i}$
   \dotfill & p.~\pageref{defgamma_i} \\
$F_i$ & $=\pi^*(E_i)$, exceptional curves of $\tau$
   \dotfill & p.~\pageref{defF_i} \\
$D$ & $=\sum_{i=1}^n (\alpha_i/2-1) E_i^*$
   \dotfill & p.\ \pageref{defD} \\
$D^*$ & $=\pi^*D=\tau^*K_{X_0}-K_X$, adjunction condition div.
   \dotfill & \fref{cd} \\
$F$ & fiber cycle of the canonical resolution
   \dotfill & p.~\pageref{defF} \\
$\bar F_i$ & exceptional curves of $\bar\tau$
   \dotfill & p.~\pageref{defbarF_i} \\
$\bar F$ & fiber cycle of the minimal resolution
   \dotfill & \fref{minfybrecycle} \\
$Z$ & fundamental cycle of the canonical resolution
   \dotfill & \fref{F} \\
$\bar Z$ & fundamental cycle of the minimal resolution
   \dotfill & \fref{barF}
\end{tabular}

\noindent
\begin{tabular}{l@{\,\,}l@{}r}
\small\it Symbol & \small\it Meaning & \small\it Ref. \\ \hline \\[-3mm]
{\it defective\,\,\,\,\,\,\,\,} & $q_i$ such that there exists $q_j>^1q_i$ with
$\alpha_j>\alpha_i\,$ $\ldots\ldots\ldots$
   \dotfill & p.\ \pageref{defdef} \\
$\Def$ & index set of defective points
   \dotfill & p.\ \pageref{defDef}
\end{tabular}

\section{Double covers of surfaces}

Throughout this paper a \emph{double cover} of a smooth irreducible complex
surface $Y$ is a finite, surjective, proper holomorphic map
$\pi:X\to Y$ of degree $2$ branched along a reduced curve,
where $X$ is a (normal) irreducible complex surface.
In the algebraic setting $\pi$ is a finite morphism of degree $2$.
For more details on covers and double covers the reader may consult \cite{bpv}.

Let $p$ be an isolated double point singularity of a surface.\newcit{p}
Since any isolated double point occurs, locally analytically, as a double cover of a smooth
surface, let us consider as our initial data a double cover $\pi_0:X_0\to Y_0$\newcit{pi0}
branched along the reduced curve $B_0$\newcit{B0}
defined in local coordinates by $f=0$, where $f$
is a square-free polynomial.\newcit{Y0}
If $x, y$ are local analytic coordinates at a point $q\in Y_0$,
then $X_0$ is defined by an equation $z^2=f(x,y)$ and $X_0$ is normal.
\newcit{X0}

If $q\notin B_0$, then $f(q)\ne0$ and there are two
pre-images of $q$ in $X_0$; at each of these points
$X_0$ is smooth and $\pi_0: X_0\to Y_0$ is unramified.
If $q\in B_0$, then there is a
single point of $X_0$ lying above $q$; $X_0$ is smooth at this point if
and only if $B_0$ is smooth at $q$.
The geometry of a smooth double cover is well-known:

\begin{lemma}
\label{doublecovers}
Let $Y$ be a smooth surface.
Let $\pi:X\to Y$ be a double cover, branched along the smooth and
reduced curve $B$ in $Y$.
\begin{enumerate}
\setlength{\parskip}{0cm}
\setlength{\itemsep}{0cm}

\item If $C$ is an irreducible component of $B$, then
$\pi_{|D}:D=\pi^{-1}(C)\to C$ is an isomorphism, $\pi^*(C)=2D$ and
$D^2=C^2/2$.

\item If $C$ is an irreducible curve in $Y$
which meets $B$ transversally in $2k$ points,
where $k \geq 1$,
then $D=\pi^{-1}(C)$ is a smooth irreducible curve on $X$ and
$\pi_{|D}:D\to C$ is a double cover branched along
the $2k$ points of intersection of $C$ with $B$.
Moreover $\pi^*(C)=D$ and $D^2=2C^2$.
\end{enumerate}
Let $E$ be a smooth rational curve in $Y$ which is not part of the branch locus
$B$. Let  $\Gamma=B_{|E}$ be the divisor of  intersection of $B$ with $E$.
\begin{enumerate}
\addtocounter{enumi}{2}
\setlength{\parskip}{0cm}
\setlength{\itemsep}{0cm}

\item If $\Gamma$ is an \emph{even} divisor, say $\Gamma=2Q$
(in particular if $\Gamma=Q=0$),
then $\pi^{-1}(E)=D_1+D_2$, where
$\pi_{|D_i}:D_i\to E$ is an isomorphism and $D_i^2=E^2-\deg(Q)$
for $i=1,2$, and $D_1\cdot D_2=\deg(Q)$.

\item If $\Gamma$ is \emph{not} even, then $D=\pi^{-1}(E)$ is
an irreducible curve, $\pi_{|D}:D\to E$ is a double cover
and $D^2=2E^2$.
$D$ is singular at those points $q \in E$
where $\Gamma(q) \geq 2$; locally near $p = \pi^{-1}(q)$
the curve $D$ has the analytic equation $z^2 = x^n$,
where $n = \Gamma(q)$.
\end{enumerate}
\end{lemma}

\begin{pf}
One can check locally the properties of $\pi$. If $C$ is not contained in $B$,
then
$\pi_{|\pi^{-1}(C)}:\pi^{-1}(C)\to C$ is a double cover branched along
$B_{|C}$ and
it is surely reducible if $B_{|C}=2Q$ for some divisor $Q\ne0$ on $C$. Since
$\pi$ has degree 2, we see that intersections double
after applying $\pi^*$, i.e.,  for any divisors $C_1$ and $C_2$ in $Y$,
$(\pi^*C_1\cdot\pi^*C_2)_X=2(C_1\cdot C_2)_Y$,
where $(\,\cdot\,)_X$ (resp.\ $(\,\cdot\,)_Y$) is the intersection form
in $X$ (resp.\ $Y$). So all the claims about $D^2$ are trivial. Regarding again
point 3, note that an
unramified double cover of $\PP^1$ is reducible, by Hurwitz formula (or
the simply-connectedness of $\PP^1$).\qed
\end{pf}

In the assumption of Lemma \ref{doublecovers},
$B$ is an even divisor and
\begin{equation}\label{piOx}
\pi_*\CO_X \cong \CO_Y \oplus \CO_Y(-B/2),
\end{equation}
while Riemann-Hurwitz Formula is:
\begin{equation}
\label{Riemann-H}
K_X=\pi^*(K_Y\otimes\CO_Y(B/2))
\end{equation}
(see \cite[Lemmas 17.1 and 17.2]{bpv}).
If $\pi:X\to Y$ is a double cover with $X$ normal but not smooth, one can
still define
the canonical divisor $K_X$ and the same formula \fref{Riemann-H} holds. For
a normal surface $X$ the canonical divisor $K_X$ is defined as the Weil
divisor class $\Div(s)$, where $s$ is a rational canonical differential.
Since Weil divisors
on a normal scheme do not depend on closed subsets of $\codim \geq 2$,
one can easily verify that
(\ref{Riemann-H}) holds
by considering the smooth double cover $\pi_{|{X_{sm}}}:X_{sm}\to Y\setminus\sing(B)$,
where $X_{sm}=X\setminus \pi^{-1}(\sing(B))$.

\section{Resolving a singular double cover}

Suppose that the branch curve $B_0$ of $\pi_0:X_0\to Y_0$ is not smooth
at the point $q_1$, thus $p=\pi_0^{-1}(q_1)$ is a double point
singularity of $X_0$.\newcit{q1}

Let us define $\mu_1=\mult_{q_1}(B)=2k+\veps_1$, with $\veps_1\in\{0,1\}$.
In order to get a resolution of the singularity $p\in X_0$,
we begin with resolving the branch locus $B_0$ of $\pi_0$.
Let $\sigma:Y_1 \to Y_0$ be the blowing-up at $q_1 \in Y_0$,
with exceptional curve $E_1$, and let $\tilde{B}_1$ be the proper
transform of $B_0$
in $Y_1$.

\begin{lemma}
\label{1blowup}
Let $X_1$ be the double cover of $Y_1$ branched along
$\tilde{B}_1 + \varepsilon_1 E_1$.
Then $X_1$ is the normalization of the pullback $X_0 \times_{Y_0} Y_1$,
and as such is both a double cover of $Y_0$
and dominates the singular surface $X_0$.
\end{lemma}

\begin{pf}
The pullback $X_0\times_{Y_0} Y_1$,
is a double cover of $Y_1$ (using the second projection)
branched along the pullback of $B_0$, which is
$\pi^*(B_0) = \mu_1 E_1+\tilde{B}_1$
(in fact it is defined by $z^2 = f(x,y)$ also).
However in the local coordinates $(u,v)$ of $Y_1$,
where $x = u$ and $y = uv$,
we have that $f(x,y) = f(u,uv) = u^{\mu_1} g(u,v)$
where $g(u,v)$ is a function
whose series expansion is not divisible
by $u$, and defines the proper transform $\tilde{B}_1$ of $B_0$.
Then the pullback is defined by $z^2 = u^{\mu_1} g(u,v)$,
and is not normal if $\mu_1 \geq 2$;
if $\mu_1 = 2k+\veps_1$ with $\veps \in \{0,1\}$,
then $w = z/u^k$  satisfies the monic equation
$w^2 = u^{\veps_1} g(u,v)$.
The normalization $X_1$ of the pullback is defined by this monic
equation,
and is clearly a double cover of $Y_1$
branched along $\tilde{B}_1 + \veps_1 E_1$.
\qed
\end{pf}

This normal surface $X_1$ dominates $X_0$, via the first projection,
and gives a partial resolution of the double point singularity.
We have passed to the double cover $X_1 \to Y_1$,
and may iterate the procedure,
continuing to blow up the branch curve at each of its singular points,
then normalizing the double cover equation. It is already known that
this process eventually terminates in a smooth double cover
$X_n\to Y_n$ that is called the \emph{canonical} resolution of
$X_0$ (see Theorem \ref{thmcanres}).
Lemma \ref{1blowup} says that the series of double covers are
determined by the parities of the multiplicities of the singular points
involved.
If the multiplicity is even, then the exceptional curve is not
part of the new branch locus, only the proper transform;
if the multiplicity is odd, then the exceptional curve is
part of the new branch locus.
This is all well-known, see for example \cite{bpv}.
We note that the multiplicity $\mu_1$ may be determined on $Y_1$
by $\mu_1 = \tilde{B}_1\cdot E$.

\section{Blowing up a smooth surface}\label{blowups}

Let us consider a sequence of blowing-ups
of a smooth surface $Y_0$, each one at a single point.
We fix a particular order for the blowing-ups,
and let $Y_i$ be the surface obtained after the $i$-th blowing-up
$\sigma_i:Y_i\to Y_{i-1}$\newcit{sigma_i}
at a point $q_i\in Y_{i-1}$.\newcit{q_i}
Let $\sigma_{ij}$ be the composition
of the blowing-up maps from $Y_j$ to $Y_i$:
\begin{equation}\label{sigma_ij}
\sigma_{ij}=\sigma_j\circ\cdots\circ\sigma_{i+2}\circ\sigma_{i+1}:Y_j\to
Y_i,
\end{equation}
for $0\le i<j \le n$, where $n$ is the total number of blowing-ups.
Set $\sigma=\sigma_{0n}:Y_n\to Y_0$, $Y=Y_n$
and $(\,\cdot\,)_i$, resp.\ $(\,\cdot\,)$\newcit{()i},
the intersection form in $Y_i$, resp.\ in $Y_n$.
The exceptional curve $E_i=\sigma_i^{-1}(q_i)$\newcit{defEi}
in $Y_i$ satisfies
$(E_i\cdot E_i)_i=-1$ and $(E_i\cdot\sigma^*_i(C))_i=0$ for any divisor $C$ of $Y_{i-1}$.
We will abuse notation
and refer to the proper transform of $E_i$ on $Y_{j\geq i}$
also as $E_i$, so $E_1,\ldots,E_n$ are the \emph{exceptional
curves} for $\sigma$.
Let $E_i^*=\sigma^*_{in}(E_i)$\newcit{defE*i}
be the total  transform
of $E_i$ in $Y_n$ via $\sigma_{in}$.

It is well-known that the relative Picard group $\Pic(Y_n)/\sigma^*\Pic(Y)$ is freely
generated by the classes $\{E_i\}_{1\le i\le n}$, as well as by
$\{E_j^*\}_{1\le j\le n}$,
and the latter ones are an \emph{orthonormal} basis in the sense that:
\begin{equation}
\label{Ei*2}
(E_i^*\cdot E_j^*)=-\delta_{ij},
\end{equation}
where $\delta$ is the Kronecker delta. Therefore we may write:
\begin{equation}
\label{Ei}
E_i=\sum_{j=1}^n n_{ij}E_j^*,\qquad E^*_j=\sum_{k=1}^n m_{jk}E_k
\end{equation}
where the matrix $M=(m_{jk})$ is the inverse of $N=(n_{ij})$.
Let $Q=(q_{ij})$ be the strictly upper triangular matrix defined by
$q_{ij}=1$ if $q_j$ lies on $E_i$ and $q_{ij}=0$ otherwise.\newcit{defQ}
Following the classical terminology, we say that the point $q_j$
is \emph{proximate} to $q_i$, and we write $q_j\prox q_i$, if and
only if $q_{ij}=1$.\newcit{qjproxqi}
\begin{lemma}\label{lemmaNQ}
Let $I$ be the identity matrix. Then
$M=I+Q+\cdots+Q^{n-1}$ and
\begin{equation}\label{N=I-Q}
N=I-Q.
\end{equation}
\end{lemma}

\begin{pf}
The first formula follows from \fref{N=I-Q}, because $Q^m=0$ for $m\ge n$.
We prove \fref{N=I-Q} by induction on $n$. For $n=1$, it is clear. If
\fref{N=I-Q} holds for $n-1$, then the proper transform of $E_i$ in $Y_{n-1}$
is
$
\sigma^*_{i,n-1}(E_i)-\sum_{j=i+1}^{n-1}q_{ij}\sigma^*_{j,n-1}(E_j),
$
so the proper transform of $E_i$ in $Y_n$ is
$
\sigma^*_n(E_i)-q_{in}E_n=\sigma^*_{in}(E_i)
-\sum_{j=i+1}^{n-1}q_{ij}\sigma^*_{jn}(E_j)-q_{in}E_n^*
$
and we conclude comparing with the first formula in  \fref{Ei}.
\qed
\end{pf}

More explicitly, $M$ can be computed inductively from $Q$ as
follows. Suppose that the first $n-1$ columns of $M$ are known. If
$q_{in}$ is the unique non-zero entry of the last column of $Q$,
then the last column of $M$ is equal to the $i$-th column of $M$
(apart $m_{nn}=1$), otherwise there is also $q_{jn}=1$ and the
last column of $M$ is the sum of the $i$-th and the $j$-th column
of $M$ (except $m_{nn}=1$).
Let us consider the matrix $Q$ as the adjacency matrix
of a directed graph $G$,
that we call the \emph{Enriques} digraph of $\sigma$:
the vertices of $G$ are the points $q_i$, for $i=1,\ldots,n$,
and there is an arrow from $q_j$ to $q_i$ if and only if $q_{ij}=1$
(i.e.\ $q_j$ is proximate to $q_i$).

\begin{remark}
The properties of $Q$ imply that an Enriques digraph is
characterized by the following four properties
(see \cite[\S5]{c}, \cite[pp.\ 213--214]{kp}):
\begin{enumerate}
\renewcommand{\theenumi}{\roman{enumi}}
\renewcommand{\labelenumi}{\theenumi)}
\setlength{\parskip}{0cm}
\setlength{\itemsep}{0cm}

\item there is no directed cycle;

\item every vertex has \emph{out-degree} at most 2;

\item if $q_i\to q_j$ and $q_i\to q_k$, with $j\ne k$,
then either $q_j\to q_k$ or $q_k\to q_j$;

\item there is at most one $q_i$ with $q_i\to q_j$
and $q_i\to q_k$, if $j\ne k$.
\end{enumerate}
\end{remark}

Note that the \emph{in-degree} of $q_i$
(the number of arrows ending in $q_i$) is $-E_i^2-1$.

Since we need only to resolve an isolated singularity at $q_1$, we
assume to blow up only points lying on the total exceptional divisor,
i.e., we assume that $q_i\in\sigma_{1,i-1}^{-1}(q_1)$, for every $i>1$.
This means that only the first column $Q_1$ of $Q$ is everywhere zero,
so the Enriques digraph is connected.
Recall that a point $q_j$ is called \emph{infinitely near} to $q_i$, and we write
$q_j>q_i$, if $q_j\in\sigma_{i,j-1}^{-1}(q_i)$.
Thus each $q_i$ is infinitely near to $q_1$.
Let us define the \emph{infinitesimal order} inductively.
If $q_j>q_i$ and there is no $q_k$ such that $q_j>q_k>q_i$,
then $q_j$ is infinitely near of \emph{order one} to $q_i$
and we write $q_j>^1q_i$.
If $q_j>^1q_k>q_i$, then by induction $q_k>^mq_i$ for some $m$
and we set $q_j>^{m+1}q_i$.\newcit{qj>qi}

Usually, the main combinatorial tool used for blowing-ups is the dual
graph of the exceptional curves and their self-intersection
numbers, that are the entries of the intersection matrix
$S=(s_{ij})$, where $s_{ij}=(E_i\cdot E_j)$.\newcit{defS}

\begin{lemma}
The configuration of the exceptional curves $E_i$ of $\sigma$ may be given by
only one of the following matrices: $M$, $N$, $Q$, or $S$.
Indeed anyone of them determines canonically all the others.
\end{lemma}

\begin{pf}
Recall that $N=M^{-1}=I-Q$. Formulas (\ref{Ei}) and (\ref{Ei*2}) imply that:
\[
s_{ij}=\left(\sum_{k=1}^nn_{ik}E_k^*\cdot\sum_{h=1}^nn_{jh}E_h^*\right)
=\sum_{k=1}^n\sum_{h=1}^nn_{ik}n_{jk}(E_k^*\cdot E_h^*)
=-\sum_{k=1}^nn_{ik}n_{jk},
\]
so $S=-NN^T=N(-I)N^T$, that is the decomposition of $S$ in an
unipotent upper triangular, a diagonal and an unipotent lower
triangular matrix. Such a decomposition is known to be unique by
linear algebra.\qed
\end{pf}


\section{The proper transform of the singular curve}\label{propertransform}

We now consider a reduced curve $B_0$ on $Y_0$
with a singular point at the point $q_1 \in Y_0$
which is being blown up.
Let us denote with $\tilde B_i$\newcit{tildeBi}
the proper transform of $B_0$ in $Y_i$
via the sequence of blowing-ups $\sigma_{0i}$, for every
$i=1,\ldots,n$. Recall that the following formula holds in $\Pic Y_n$:
\begin{equation}
\label{tildealfaj}
\sigma^*(B_0)=\tilde B_n+\sum_{i=1}^n\tilde\alpha_iE_i^*,\qquad
\hbox{where $\tilde\alpha_i=\mult_{q_i}(\tilde B_{i-1}).$}
\end{equation}
Abusing language a little, usually $\tilde\alpha_i$ is called the
\emph{multiplicity} of $B_0$ at $q_i$.
Note that $\tilde\alpha_i$ may be determined also in $Y_i$ by
$\tilde\alpha_i=(\tilde B_i\cdot E_i)_i$.

On the other hand, in $\Pic Y_n$ we may write also:
\begin{equation}
\label{tildebetaj}
\sigma^*(B_0)=\tilde B_n+\sum_{i=1}^n\tilde\beta_iE_i.
\end{equation}
for some non-negative integers $\tilde{\beta}_i$.
Putting the second formula of (\ref{Ei}) in (\ref{tildealfaj})
we find that the $\tilde\beta_i$'s can be computed from the $\tilde\alpha_i$'s as
follows:
\begin{equation}
\label{beta=alfaM}
\tilde\beta_i=\sum_{j=1}^n\tilde\alpha_jm_{ji},\qquad\hbox{or shortly}\qquad
\tilde\beta=\tilde\alpha M,
\end{equation}
where $\tilde\alpha$ and $\tilde\beta$ are row vectors with the obvious entries.

The quantities $\tilde\alpha$ and $\tilde\beta$ may also be
determined on $Y_n$ knowing how $\tilde B_n$ intersects the
exceptional curves $E_i$. Indeed, intersecting (\ref{tildebetaj})
with $E_i$ gives
$
0=\tilde B_n\cdot E_i+\sum_{j=1}^n\tilde\beta_jE_i\cdot E_j
$\,
that is, setting $\tilde\gamma_j=(\tilde B_n\cdot E_j)$\newcit{tildegammaj}
and $\tilde\gamma$ the corresponding
row vector:
\begin{equation}
\label{gamma=alfaM}
\tilde\gamma=-\tilde\beta S=\tilde\beta NN^T=\tilde\alpha N^T,
\qquad\hbox{so}\qquad\tilde\alpha=\tilde\gamma M^T.
\end{equation}
Note that $\tilde B_n$ satisfies $(\tilde B_n\cdot E_i)\ge0$ for every $i$,
which is equivalent by \fref{tildealfaj},  \fref{Ei}, \fref{N=I-Q} and
(\ref{Ei*2}) to the so-called
\emph{proximity inequality} at $q_i$:
\begin{equation}
\label{disprox1}
\tilde\alpha_i\ge\sum_{j=1}^nq_{ij}\tilde\alpha_j=\sum_{j:q_j\prox
q_i}\tilde\alpha_j.
\end{equation}

Suppose that $B_0$ has only an isolated singularity in $q_1$ and that
$\sigma:Y_n\to
Y_0$ resolves the singularities of $B_0$, i.e., the proper transform $\tilde
B_n$ is
smooth. Then the topological type of the singularity of $B_0$ at $q_1$ is
completely
determined by the matrix $M$, which carries the configuration of the exceptional
curves $E_i$ for $\sigma$, and
the intersection divisor $\tilde{\Gamma}_i=\tilde B_{n|E_i}$,\newcit{tildeGammaj}
which says
how $\tilde B_n$ meets $E_i$, for every $i=1,\ldots,n$. Each $\tilde{\Gamma}_i$ is a
non-negative divisor on $E_i$ and if $q$ is a point of intersection of two
exceptional
curves $E_i$ and $E_j$, then:
\begin{equation}
\label{B|EiEj}
\tilde{\Gamma}_i(q)=0\iff\tilde \Gamma_j(q)=0.
\end{equation}
Moreover if these numbers are non-zero, then at least one number
is equal to one. Condition (\ref{B|EiEj}) says that
either $\tilde B_n$ passes through $q$ or not, while the latter
statement that $\tilde B_n$ cannot be tangent in $q$ to both
exceptional curves simultaneously.
These divisors $\tilde \Gamma_i$ express the combinatorial information
of the singularity of the branch curve completely. Given the configuration
of the exceptional curves, they
can be arbitrary, subject to the above condition.

We remark that the knowledge of $M$ and of the degrees
$\tilde\gamma_i=\deg(\tilde\Gamma_i)$ is equivalent to the knowledge
of $Q$ and the multiplicities $\tilde\alpha_i$ of $B_0$ at $q_i$, by \fref{gamma=alfaM}.

Therefore we will define the \emph{weighted Enriques digraph}
of $q_1\in B_0$ by attaching to each vertex $q_i$ of the Enriques
digraph the weight $\tilde\alpha_i=\mult_{q_i}(\tilde B_{i-1})$.

\begin{example}
\label{example1}
Let $B_0\subset Y_0$ be defined locally near the origin $q_1=(0,0)$ by:
\[
x(y^2-x)(y^2+x)(y^2-x^3)(y^2+x^3)=0.
\]
\end{example}

Clearly $B_0$ has multiplicity $\tilde\alpha_1=7$ at $q_1$.
Blow up $q_1$: $\tilde B_1\subset Y_1$ has two singular points $q_2$ and $q_3$ on $E_1$ of
multiplicity $\tilde\alpha_2=3$ and $\tilde\alpha_3=2$ ($q_2$ and  $q_3$
are infinitely near points of order one to $q_1$).
Then blow up $q_2$ and $q_3$. $\tilde B_2\subset Y_2$ meets transversally
$E_2$ and it is smooth at those points, but $\tilde B_3\subset Y_3$ has
in $q_4=E_1\cap E_3$ a point of multiplicity $\tilde\alpha_4=2$,
so $q_4$ is proximate to both $q_3$ and $q_1$, but $q_4$ is
infinitely near of order 2 to $q_1$. Classically, $q_4$ is called
a \emph{satellite} point to $q_1$. Finally $\tilde B_4\subset Y_4$ is smooth.

The configuration of the exceptional curves is determined by
anyone of the following matrices:
\[
Q=\begin{pmatrix}
0 & 1 & 1 & 1\\
0 & 0 & 0 & 0\\
0 & 0 & 0 & 1\\
0 & 0 & 0 & 0
\end{pmatrix}\quad
M=\begin{pmatrix}
1 & 1 & 1 & 2\\
0 & 1 & 0 & 0\\
0 & 0 & 1 & 1\\
0 & 0 & 0 & 1
\end{pmatrix}\quad
S=\begin{pmatrix}
-4 & 1 & 0 & 1\\
1 & -1 & 0 & 0\\
0 & 0 & -2 & 1\\
1 & 0 & 1 & -1
\end{pmatrix}
\]
and the combinatorial data of $\tilde B_4$ by anyone of the
following vectors:
\[
\tilde\alpha=(7,3,2,2),\qquad\tilde\beta=(7,10,9,18),\qquad
\tilde\gamma=(0,3,0,2),
\]
that we encode in the following weighted Enriques digraph:
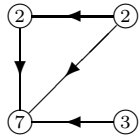
\[
\begin{picture}(60,60)
\thinlines
\put(14.5,10){\line(1,0){31}}
\put(10,14.5){\line(0,1){31}}
\put(14.5,50){\line(1,0){31}}
\put(13.5,13.5){\line(1,1){33}}
\put(10,50){\circle{9}}
\put(50,10){\circle{9}}
\put(10,10){\circle{9}}
\put(50,50){\circle{9}}
\put(48,7.8){$\scriptstyle3$}
\put(8,7.8){$\scriptstyle7$}
\put(8,47.8){$\scriptstyle2$}
\put(48,47.8){$\scriptstyle2$}
\thicklines
\put(27,10){\vector(-1,0){0}}
\put(10,27){\vector(0,-1){0}}
\put(27,50){\vector(-1,0){0}}
\put(27,27){\vector(-1,-1){0}}
\end{picture}
\]

\section{Resolving the branch locus of a double cover}\label{resolving}

We return to consider a normal double cover $\pi_0:X_0\to Y_0$ branched along
the singular reduced curve $B_0$.
We desingularize $B_0$ by successive blowing-ups as in the previous sections.
We have seen that $\sigma$ induces a normal double cover $\pi_n:X_n\to Y_n$, which is the
normalization of the pullback $X_0\times_{Y_0}Y_n$.
As in Lemma \ref{1blowup}, $\pi_n$ is branched along a
reduced curve $B_n$, obtained from $\sigma^*(B_0)$ by removing the
multiple components an even number of times.
Hence $B_n$ is made of the proper transform $\tilde B_n$ of $B_0$ and possibly
of some exceptional curves $E_i$.
We set $\veps_i$ equal to one or zero, depending on whether $E_i$
is part of the branch locus of $Y_n$
or not (and we say that $E_i$ is \emph{branched}, resp.\ \emph{unbranched}).
Setting $\veps$ the corresponding row
vector, by (\ref{tildebetaj}) and (\ref{beta=alfaM}) we have that:
\begin{equation}
\label{p=beta}
\veps=\tilde\beta\bmod2\quad\hbox{and}\quad\veps=\tilde\alpha M\bmod2.
\end{equation}
Note that the branchedness of $E_i$ is
determined  at the moment $E_i$ is created on $Y_i$ by blowing up the point
$q_{i}\in Y_{i-1}$. Indeed, Lemma \ref{1blowup} says that
$\veps_i=1$
(resp.~$\veps_i=0$) if the multiplicity $\mu_i$ at $q_i$ of the branch locus $B_{i-1}$ of
$\pi_{i-1}:X_{i-1}\to Y_{i-1}$ is odd (resp.~even).\newcit{pi_i}
Shortly:
\begin{equation}
\label{p=alfa'}
\veps=\mu\bmod2.
\end{equation}
Assuming inductively to know
$\veps_j$
for $j<i$, the multiplicity of $B_{i-1}$ at $q_i$ is
\begin{equation}
\label{alfa'}
\mu_i=\tilde\alpha_i+\sum_{j=1}^{i-1}\veps_iq_{ji}=
\tilde\alpha_i+\sum_{j:q_i\to q_j}\veps_j,\qquad\hbox{or shortly}\qquad
\mu=\tilde\alpha +\veps Q.
\end{equation}

If we have blown up to make the proper transform $\tilde{B}_n$ smooth,
we still may not have the total branch locus $B_n$ smooth.
Singularities of the total branch locus now come from two
sources: intersections of $\tilde{B}_n$ with branched $E_i$'s,
and intersections between two different branched $E_i$'s.
We first take up the former case.
Suppose that $\tilde{B}_n$ meets the exceptional configuration $\bigcup_i E_i$
at a point $q_{n+1}$.
We blow up $q_{n+1}$ to create a new surface $Y_{n+1}$,
a new proper transform $\tilde{B}_{n+1}$,
and a new exceptional curve $E_{n+1}$.
Since $\tilde{B}_n$ is smooth at $q_{n+1}$,
we have $(\tilde{B}_{n+1} \cdot E_{n+1})_{n+1} = 1$.
We now have new intersection divisors $\tilde\Gamma_i'$
for each $i = 1,\ldots,n+1$.
These are related to the previous intersection divisors $\tilde {\Gamma}_i$'s
for $i=1,\ldots,n+1$ as follows:
\[
\tilde{\Gamma}'_i=
\begin{cases}
\tilde\Gamma_i & \text{if $q_{n+1}\notin E_i$}\\
\tilde{\Gamma}_i-q_{n+1} & \text{if $q_{n+1}\in E_i$}
\end{cases}
\]
for every $i=1,\ldots,n$ and $\tilde{\Gamma}'_{n+1}=q$, where $q$ is the point of
intersection of $\tilde B_{n+1}$ with $E_{n+1}$.
We may now iterate this construction,
and arrive at the situation
(increasing the number of blowing-ups $n$)
that the proper transform $\tilde{B}_n$
does not meet any branched exceptional curves, i.e., the new exceptional
divisors $\tilde{\Gamma}'_i$ are zero for each $i$ such that $\veps_i=1$.
Finally if any two exceptional curves now meet,
we simply blow up the point of intersection once,
and obtain an unbranched exceptional curve
which now separates the two branched exceptional curves.

At this point we have a nonsingular total branch locus, hence a smooth  double cover.
We note that we still have the matrices $M, N, Q, S$ for the current
configuration of exceptional curves
and the numbers $\tilde\alpha_i$, $\tilde\beta_i$, $\tilde\gamma_i$ and $\veps_i$,
defined for each $i$, as before.
All the above process gives a proof of the following theorem
(cf.\ \cite[III.\S6]{bpv}, \cite[Theorem 3.1]{laufer3}):

\begin{theorem}[The canonical resolution]
\label{thmcanres}
Let $\pi_0:X_0\to Y_0$ be a double cover with $X_0$ normal and $Y_0$ smooth.
Then
there exists a birational morphism $\sigma:Y\to Y_0$ such that the normalization
$X$
of the pullback $X_0\times_{Y_0}Y$ is smooth. Moreover $\pi:X\to Y$ is a double
cover and the diagram
\begin{equation}
\label{canonicaldiagram}
\begin{array}{ccc}
X & \mapright\tau & X_0 \\ \mapdown\pi & & \mapdown{\pi_0} \\ Y &
\mapright\sigma &
Y_0
\end{array}
\end{equation}
commutes.
So $\tau:X\to X_0$ is a resolution of singularities of $X_0$.\qed
\end{theorem}
We say that $\tau:X\to X_0$\newcit{tau}
(toghether with the double cover map $\pi:X\to Y$\newcit{pi:})
is the \emph{canonical} resolution of the double
cover $\pi_0:X_0\to Y_0$, because $X$ and $Y$ are unique, up to isomorphism, assuming
that
the centers of the blowing-ups $\sigma_i:Y_i\to Y_{i-1}$, which
factorizes $\sigma$, are always singular points of the branch curve
of $X_{i-1}\to Y_{i-1}$.
We will see in section \ref{minimalres} that the
canonical resolution might not be minimal. However, we may think the canonical
resolution as the ``minimal'' resolution in the category of double covers
over smooth surfaces.

\section{The branch curve of the canonical resolution}\label{branchcurve}

As in the previous section, let $\tau:X\to X_0$ be the canonical resolution of
$X_0$
and $\pi:X\to Y$ be the smooth double cover. Let us write $B$ for the branch
curve of
$\pi:X\to Y$ and suppose that $\sigma:Y\to Y_0$ factorizes in $n$ blowing-ups
 $\sigma_i:Y_i\to Y_{i-1}$, with $Y_n=Y$,\newcit{Y}
so the exceptional curves for
$\sigma$ are $E_1,\ldots,E_n$ and all the formulas in the previous sections hold for
$B_n=B$\newcit{defB}
and $\tilde B_n=\tilde B$. The branch curve $B$ can be written in $\Pic Y$ as:
\begin{equation}
\label{Balfabeta} B=\tilde B+\sum_{i=1}^n\veps_iE_i
=\sigma^*(B_0) - \sum_{i=1}^n\beta_iE_i
=\sigma^*(B_0) - \sum_{i=1}^n\alpha_iE_i^*
\end{equation}
for some non negative integers $\beta_i$ and $\alpha_i$.
Comparing (\ref{Balfabeta}) with (\ref{tildebetaj}),
(\ref{tildealfaj}) and \fref{Ei}
we see that:
\begin{equation}
\label{alfa=alfa'-p}
\beta=\tilde\beta-\veps
\qquad\hbox{and}\qquad
\alpha=\tilde\alpha+\veps N=
\tilde\alpha +\veps (Q-I).
\end{equation}
Moreover, from formulas (\ref{alfa=alfa'-p}), (\ref{beta=alfaM})
and (\ref{alfa'}) we find that:
\begin{equation}
\label{alfabetamu}
\alpha=\mu-\veps  \qquad\hbox{and}\qquad
\alpha=\beta N.
\end{equation}
By \fref{alfa=alfa'-p} and (\ref{p=beta}),  \fref{alfabetamu} and \fref{p=alfa'} it
follows that  the $\beta_i$'s and the  $\alpha_i$'s are all
even. From (\ref{alfa'}) and (\ref{alfa=alfa'-p})
it is clear that
$
\alpha=\mu=\tilde\alpha$ if and only if $\veps=0,
$
that happens if $\tilde\alpha_i$ is even for every $i=1,\ldots,n$,
i.e.\ if $B_0$ has even multiplicity at each singular point,
including the infinitely near ones.

In order to measure how the branch curve $B$ of the canonical
resolution meets the exceptional curves, let us introduce
the following intersection divisors:
\begin{equation}\label{Wi}
\Gamma_i=B_{|E_i}=\tilde B_{|E_i}+\sum_{j:\veps_j=1}{E_j}_{|E_i}
\end{equation}
and set $\gamma_i=\deg(\Gamma_i)$.\newcit{defgamma_i}
By definition
$\Gamma_i=0$ if $E_i$ is branched. However, $\Gamma_i$ could be zero even if
$\veps_i=0$.
We claim that the $\Gamma_i$'s have the following properties:
\begin{enumerate}
\setlength{\parskip}{0cm}
\setlength{\itemsep}{0cm}

\item $\Gamma_i$ is a non-negative divisor and $\gamma_i$ is even;

\item if $q=E_i\cap E_j$ and $\veps_i=\veps_j=0$, then $\Gamma_i(q)=0$
if and only if $\Gamma_j(q)=0$.
If these numbers are non-zero, then at least one of them is
equal to 1.

\item if $q=E_i\cap E_j$, $\veps_i=1$ and $\veps_j=0$, then $\Gamma_j(q)=1$.
\end{enumerate}

It suffices to show that $\gamma_i$ is even, since all the other
properties of $\Gamma_i$ are induced by those of $\tilde{\Gamma}_i$,
which we have already seen in section \ref{propertransform}.
If $E_i$ is branched, then $\gamma_i=0$ and the thesis is trivial.
Otherwise, if $E_i$ is unbranched, then:
\[
\gamma_i=\tilde\gamma_i+\sum_{j\ne i}\veps_j(E_i\cdot E_j)\equiv
\tilde\gamma_i+\sum_{j=1}^n\tilde\beta_js_{ij}\bmod 2
\]
hence $\gamma\equiv\tilde\gamma+\tilde\beta S\bmod2$ and the claim follows from
$\tilde\gamma=-\tilde\beta S$ (see (\ref{gamma=alfaM})).

In order to encode the combinatorial data of the double point
singularity $p=\pi_0^{-1}(q_1)\in X_0$,
we define the \emph{weighted Enriques digraph} of $p$
by attaching to each vertex $q_i$ of the Enriques digraph
of $q_1\in B_0$ the weight $\mu_i=\mult_{q_i}(B_{i-1})$.

In the next example we will illustrate in detail how the canonical
resolution process goes on.

\begin{example}
\label{example3}
Let $B_0$ be a curve defined locally near the origin $q_1=(0,0)$ by:
\[
y(y^2-x^3)=0.
\]
\end{example}

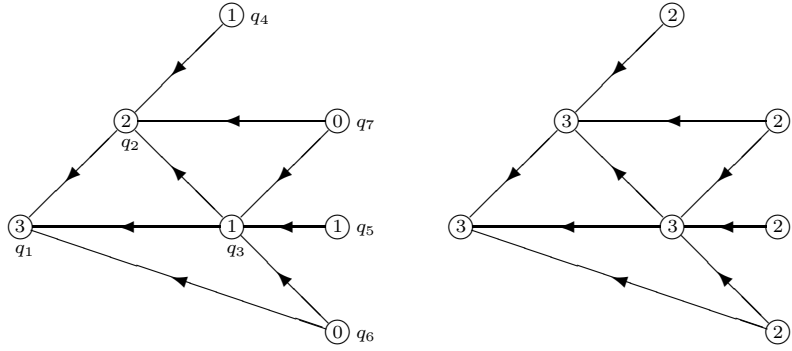
\begin{figure}[ht]
\centerline{
\begin{picture}(140,140)
\thinlines
\put(10,50){\circle{9}}
\put(50,90){\circle{9}}
\put(90,50){\circle{9}}
\put(90,130){\circle{9}}
\put(130,10){\circle{9}}
\put(130,50){\circle{9}}
\put(130,90){\circle{9}}
\put(14.5,50){\line(1,0){71}}
\put(94.5,50){\line(1,0){31}}
\put(54.5,90){\line(1,0){71}}
\put(13.5,53.5){\line(1,1){33}}
\put(53.5,93.5){\line(1,1){33}}
\put(53.5,86.5){\line(1,-1){33}}
\put(93.5,53.5){\line(1,1){33}}
\put(93.5,46.5){\line(1,-1){33}}
\put(14.2,48.6){\line(3,-1){111.6}}
\put(8,47.8){$\scriptstyle3$}
\put(48,87.8){$\scriptstyle2$}
\put(88.1,47.8){$\scriptstyle1$}
\put(88.1,127.8){$\scriptstyle1$}
\put(128.1,7.8){$\scriptstyle0$}
\put(128.1,47.8){$\scriptstyle1$}
\put(128.1,87.8){$\scriptstyle0$}
\put(8,40){$\scriptstyle q_1$}
\put(48,80){$\scriptstyle q_2$}
\put(88,40){$\scriptstyle q_3$}
\put(97,127.8){$\scriptstyle q_4$}
\put(137,7.8){$\scriptstyle q_6$}
\put(137,47.8){$\scriptstyle q_5$}
\put(137,87.8){$\scriptstyle q_7$}
\thicklines
\put(47,50){\vector(-1,0){0}}
\put(107,50){\vector(-1,0){0}}
\put(87,90){\vector(-1,0){0}}
\put(27,67){\vector(-1,-1){0}}
\put(67,107){\vector(-1,-1){0}}
\put(107,67){\vector(-1,-1){0}}
\put(67,73){\vector(-1,1){0}}
\put(107,33){\vector(-1,1){0}}
\put(67,31){\vector(-3,1){0}}
\end{picture}
$\qquad$
\begin{picture}(140,140)
\thinlines
\put(10,50){\circle{9}}
\put(50,90){\circle{9}}
\put(90,50){\circle{9}}
\put(90,130){\circle{9}}
\put(130,10){\circle{9}}
\put(130,50){\circle{9}}
\put(130,90){\circle{9}}
\put(14.5,50){\line(1,0){71}}
\put(94.5,50){\line(1,0){31}}
\put(54.5,90){\line(1,0){71}}
\put(13.5,53.5){\line(1,1){33}}
\put(53.5,93.5){\line(1,1){33}}
\put(53.5,86.5){\line(1,-1){33}}
\put(93.5,53.5){\line(1,1){33}}
\put(93.5,46.5){\line(1,-1){33}}
\put(14.2,48.6){\line(3,-1){111.6}}
\put(8,47.8){$\scriptstyle3$}
\put(48,87.8){$\scriptstyle3$}
\put(88.1,47.8){$\scriptstyle3$}
\put(88.1,127.8){$\scriptstyle2$}
\put(128,7.8){$\scriptstyle2$}
\put(128,47.8){$\scriptstyle2$}
\put(128,87.8){$\scriptstyle2$}
\thicklines
\put(47,50){\vector(-1,0){0}}
\put(107,50){\vector(-1,0){0}}
\put(87,90){\vector(-1,0){0}}
\put(27,67){\vector(-1,-1){0}}
\put(67,107){\vector(-1,-1){0}}
\put(107,67){\vector(-1,-1){0}}
\put(67,73){\vector(-1,1){0}}
\put(107,33){\vector(-1,1){0}}
\put(67,31){\vector(-3,1){0}}
\end{picture}}
\caption{The Enriques digraph weighted resp.\ with
the $\tilde\alpha_i$'s and the $\mu_i$'s}
\label{Fex3}
\end{figure}

Clearly, we need just two blowing-ups to smooth the proper transform
$\tilde B_2$ of $B_0$.
However, we have to blow up 5 more times in order to get a smooth
double cover.
In fact $E_1, E_2$ are branched and $\tilde B_2$ passes
through $q_3=E_1\cap E_2$ and meets $E_2$ also in another point $q_4$.
Thus $B_2=\tilde B_2+E_1+E_2$ has multiplicity $\mu_3=3$ at $q_3$,
hence $\veps_3=1$ and $\tilde B_3$ meets $E_3$ in a point $q_5$.
Now $B_3$ has only nodes in $q_4$, $q_5$, $q_6=E_1\cap E_3$
and $q_7=E_2\cap E_3$,
therefore $E_4,\ldots,E_7$ are unbranched and $B_7$ is smooth.
Our combinatorial data  are:
\[
M=\begin{pmatrix}
1 & 1 & 2 & 1 & 2 & 3 & 3\\
0 & 1 & 1 & 1 & 1 & 1 & 2\\
0 & 0 & 1 & 0 & 1 & 1 & 1\\
0 & 0 & 0 & 1 & 0 & 0 & 0\\
0 & 0 & 0 & 0 & 1 & 0 & 0\\
0 & 0 & 0 & 0 & 0 & 1 & 0\\
0 & 0 & 0 & 0 & 0 & 0 & 1
\end{pmatrix},\qquad
\begin{array}{l}
\tilde\alpha=(3,2,1,0,0,1,1)\\
\tilde\gamma=(0,0,0,1,1,0,0),\\
\mu=(3,3,3,2,2,2,2),\\
\veps=(1,1,1,0,0,0,0),\\
\alpha=(2,2,2,2,2,2,2),\\
\gamma=(0,0,0,2,2,2,2),
\end{array}
\]
that we encode in the Enriques digraph, weighted with the $\mu_i$'s
(see the right-hand side graph of Figure \ref{Fex3}).
For the readers' convenience, we inserted in Figure \ref{Fex3}
(on the left-hand side) also the Enriques digraph weighted
with the $\tilde\alpha_i$'s and labelled with the $q_i$'s.

Figure \ref{Fex3} may help to understand formula \fref{alfa'},
namely how to compute inductively the $\mu_i$'s (thus the
$\veps_i$'s) from the $\tilde\alpha_i$'s. Start from $q_1$: the weight
of $q_1$ is 3, thus $\mu_1=3$, $\veps_1=1$ and add 1 to all the
weights attached to vertices with arrows ending in $q_1$ (namely
$q_2$, $q_3$ and $q_6$). Now consider $q_2$: the actual weight of
$q_2$ is 3, hence $\mu_2=3$, $\veps_2=1$ and add 1 to the weights
of $q_3$, $q_4$ and $q_7$, which are the vertices with arrows
ending in $q_2$. Then go on inductively, for all the $q_i$'s.
Clearly, no change of weights is made at step $i$ if $\veps_i=0$.

\section{The description of the canonical resolution}\label{secdescanres}

The description of the canonical resolution of the singularity
on $X_0$ is now a combinatorial problem, using the information
of the configuration of the exceptional curves $E_i$
(described by the matrix $M$ or the Enriques digraph)
and the divisors $\tilde\Gamma_i$
(subject to the conditions stated in section \ref{propertransform}).

We then can compute the quantities $\tilde\gamma$, $\tilde\alpha$,
$\mu$, $\veps$, and determine (from $\veps$) which of the $E_i$'s
are branched curves, and finally determine the divisors $\Gamma_i$.

We now apply Lemma \ref{doublecovers} for each exceptional
curve.  If $\pi:X \to Y$ is the double cover map,
let us define $F_i=\pi^{-1}(E_i)$\newcit{defF_i}
for each $i$,
thus $F_1,\ldots,F_n$ are all the exceptional curves
for the canonical resolution $\tau:X\to X_0$.

\begin{remark}\label{remFi}
A curve $F_i$ is reducible if and only if $\veps_i=0$ and $\Gamma_i=B_{|E_i}$ is
an even divisor.
In that case, $F_i$ splits in two smooth rational curves
$F'_i$ and $F''_i$, with $F'_i\cdot F''_i=\gamma_i/2$
and $F_i^{\prime2}=F_i^{\prime\prime2}=E_i^2-\gamma_i/2$.
\end{remark}

If $\veps_i=0$, then $\pi^*(E_i)=F_i$, otherwise, if $\veps_i=1$,
then $\pi^*(E_i)=2F_i$, i.e.
\begin{equation}\label{pi*Ei}
\pi^*(E_i)=(1+\veps_i)F_i.
\end{equation}
Moreover since intersections double after applying $\pi^*$ we have that:
\begin{equation}
\label{isquare-ij}
F_i^2 = \frac{2}{{(1+\veps_i)}^2}E_i^2\qquad\hbox{and}\qquad
(F_i\cdot F_j) = (2-\veps_i-\veps_j)(E_i\cdot E_j).
\end{equation}

We claim that the arithmetic genus of $F_i$ is, for each $i$:
\begin{equation}
\label{igenus}
p_a(F_i)=\frac{\gamma_i}{2}+\veps_i-1.
\end{equation}

If $\veps_i=1$, then $F_i$
is a smooth rational curve, $\Gamma_i=0$ and the claim is trivial. If
$F_i$ splits, \fref{igenus} follows from Remark \ref{remFi}.
Otherwise, $F_i$ is a double cover of $E_i$ branched along
$\Gamma_i$ and \fref{igenus} is just Hurwitz formula.

Moreover $F_i$ is singular at a point $P$ if and only
if $\Gamma_i(\pi(P))>1$, thus in particular $F_i$ is smooth
at the intersection points with $F_j$, for each $j\ne i$.

Now we want to find an explicit formula
for the canonical divisor $K_X$.
By Riemann-Hurwitz formula (\ref{Riemann-H}) we know
that $K_X=\pi^*(K_Y+B/2)$, where
$K_Y=\sigma^*(K_{Y_0})+\sum_iE_i^*$. Therefore
by (\ref{Balfabeta}):
\begin{equation}
\label{as}
K_Y+\frac{B}{2}=\sigma^*\left(K_{Y_0}+\frac{B_0}{2}\right)-
\sum_{i=1}^n\left(\frac{\alpha_i}{2}-1\right)E_i^*.
\end{equation}
We define $D=\sum_{i=1}^n (\alpha_i/2-1) E_i^*$, so:\newcit{defD}
\begin{equation}\label{cd}
K_{X}=(\sigma \circ\pi)^*(K_{Y_0}+B_0/2)-\pi^*D=\tau^*K_{X_0}-D^*,
\end{equation}
where $D^*=\pi^*D$ is called the \emph{adjunction condition} divisor
(cf.\ section \ref{adjoints}).
We remark that $D^*\ge0$, because $\alpha_i\ge2$ for each $i$
(since in the canonical resolution process
we blow up only singular points of the branch curve).

Let us show now the explicit formula for the \emph{fiber cycle} $F$ of the
canonical resolution, already written without proof
by Franchetta and Dixon.
The \emph{fiber cycle} of $\tau$ is the maximal effective divisor $F$\newcit{defF}
contained in the \emph{scheme theoretic fiber} of $\tau$,
i.e., the subscheme of
$X$ defined by the inverse image ideal sheaf $\tau^{-1} m_{p,X_0}$ of the maximal ideal
$m_{p,X_0}$
of $p$ in
$X_0$. Therefore:
\[
F=\gcd{\{\Div(g)\,| \,\text{for all}\, g\in\tau^{-1} m_{p,{X_0}}  \}}.
\]
For example, the fiber cycle of the sequence $\sigma:Y\to Y_0$
of blowing-ups is $E_1^*$.

\begin{theorem}[Franchetta]
\label{fibercycle}
The fiber cycle of the canonical resolution is:
\begin{equation}
\label{fc}
F=\pi^*(E_1^*)=\sum_{i=1}^n m_{1i}(1+\veps_i)F_i.
\end{equation}
\end{theorem}

\begin{pf}
The second equality in \fref{fc} follows from \fref{pi*Ei} and from
the definition of $M$ (see \fref{Ei}), so it suffices to show the
first equality.
Recall that locally $x,y$ are the coordinates of $Y_0$ near $q_1=(0,0)$,
$X_0$ is defined by $z^2=f(x,y)$ and
$\pi_0$ is the projection $(x,y,z)\mapsto(x,y)$.
Hence $x$, $y$ and $z$ are
generators of the ideal sheaf $\tau^{-1}m_{p,X_{0}}$ as
$\CO_{{X_0},p}$-module.
Therefore $F$ is the greatest effective divisor contained
in the pullback divisors $\tau^*\Div x$, $\tau^*\Div y$ and
$\tau^*\Div z$. From the commutativity of the
diagram \fref{canonicaldiagram} we see that
$\tau^*\Div x={(\pi_0\circ\tau)}^*\Div x
={(\sigma\circ \pi)}^*\Div x$ and
$\tau^*\Div y={(\sigma\circ\pi)}^*\Div y$,
thus the gcd of the divisors $\tau^*\Div x$ and $\tau^*\Div y$ is the
pullback $\pi^*(E_1^*)$ of the fiber cycle $E_1^*$ of $\sigma$.
Moreover
$
\tau^*\Div z^2={(\pi_0\circ\tau)}^*\Div f={(\sigma\circ \pi)}^*\Div f,
$
thus the divisor $2\tau^*\Div z=\tau^*\Div z^2$ is equal to the pullback
$\pi^*(\sigma^*B_0)$ of the total transform $\sigma^*B_0$
which contains $E_1^*$ with multiplicity $\tilde\alpha_1\ge 2$.
So $\tau^*\Div z \supseteq \pi^*(E_1^*)$, hence
$F=\pi^*(E_1^*)$.
\qed
\end{pf}

The self-intersection of the fiber cycle is:
\begin{equation}
\label{fcsquare}
F^2 =\pi^*(E_1^*)\cdot\pi^*(E_1^*)= 2E_1^{*2}=-2=-\mult_p(X_0).
\end{equation}
By \fref{fc} and \fref{pi*Ei}, the fiber cycle $F$ has
the following properties:
\begin{equation}
\label{F*.Fi}
F\cdot F_1=(2-\veps_1)E_1^*\cdot E_1=-2+\veps_1,
\quad
F\cdot F_i=(2-\veps_i)E_1^*\cdot E_i=0
\end{equation}
for every $i>1$.
Thus the normal sheaf $\CO_{F_i}(F_i)$ of $F_i$, for $i>1$ (that can be
useful if $p_a(F_i)>0$), is given by how $F_i$ meets the other
components of $F$:
\begin{equation}
\label{Fi|Fi}
m_{1i}(1+\veps_i){F_i}_{|F_i}=-{F^{\hat\imath}}_{|F_i}
\end{equation}
where $F^{\hat\imath}=F-m_{1i}(1+\veps_i)F_i$. Finally, by
\fref{cd}, \fref{fc}
and \fref{fcsquare}, the arithmetic genus of the fiber cycle is:
\begin{equation}
\label{fibcyclegenus}
p_a(F)=(F\cdot K_X+{F}^2)/{2}+1
=E_1^*\cdot (K_Y+B/2)=\alpha_1/2-1.
\end{equation}

\section{The minimal resolution}\label{minimalres}

It may happen that the canonical resolution $\tau:X\to X_0$ of $p\in X_0$
is not minimal. However the following theorem (cf.\ \cite[Th.\ 5.4]{laufer3})
shows that the canonical resolution is not too far
to be minimal.

\begin{theorem}\label{minresthm}
Let $\tau:X\to X_0$ be the canonical resolution and
$\bar\tau:\bar X\to X_0$ the minimal one.
Then $\tau=\tau'\circ\bar\tau$, where $\tau':X\to\bar X$
is the blowing-up at finitely many distinct points.
Moreover none of these points is singular for the exceptional curves
of $\bar\tau$ and neither lies on the intersection
(necessarily transverse) of more than three of them.
\end{theorem}

The following lemma characterizes the $(-1)$-curves in $\tau^{-1}(p)$
and allows us to give an elementary proof of Theorem \ref{minresthm}, easier
than Laufer's original one.

\begin{lemma}\label{Lem10.2}
An exceptional curve $F_j$ for $\tau$ is a $(-1)$-curve if and
only if  $F_j=\pi^{-1}(E_j)$, where $E_j$ is branched and
$E_j^2=-2$. Moreover the $(-1)$-curves in $\tau^{-1}(p)$ are
disjoint.
\end{lemma}

\begin{pf}
Let $E_j$ be unbranched. By Lemma \ref{doublecovers},
if $\Gamma_j=B_{|E_j}$ is not even, then $F_j$ cannot be a $(-1)$-curve,
because $F_j^2=2E_j^2$ would be even.
If $\Gamma_j$ is even, then $\pi^{-1}(E_j)=F'_j+F''_j$
and $F'_j=F''_j=E_j^2-\deg(\Gamma_j/2)$,
thus $F'_j$ (and $F''_j$) could be a $(-1)$-curve only if $E_j^2=-1$
and $\Gamma_j=0$, that means that $q_j$ was unnecessarily blown up.
Hence we are left only with the possibility that $E_j$ is branched and
$E_j^2=2F_j^2=-2$.
Since no branched exceptional curves meet in $Y$,
neither do two $(-1)$-curves in $X$.
\qed
\end{pf}

\begin{remark}\label{rem-1}
If $\smash{F_j=\pi^{-1}(E_j)}$ is a $(-1)$-curve in $X$, then
$\smash{\mu_j}$ is odd and there is exactly one $q_i$ such
that $q_i>^1q_j$ and $\mu_i=\mu_j+1$. Usually one says
that $B_{j-1}$ has two infinitely near points of the same,
odd, multiplicity $\mu_j$ at $q_j$.
\end{remark}

\begin{pf}
Clearly $\mu_j$ is odd (and $>2$), because $\veps_j=1$. Since
$E_j$ is a $(-2)$-curve, we blew up only one point $q_i$ on
$E_j$: this means that all the intersections of $\smash{\tilde
B_j}$ with $E_j$ are supported on $q_i$, i.e.\
$\tilde\alpha_i=\tilde\alpha_j$, and the thesis follows from
\fref{alfa'}. \qed
\end{pf}

Let us denote by $\bar\tau:\bar X\to X_0$ the minimal resolution
of $p\in X_0$.

\medskip

\noindent\textbf{Proof of Theorem \ref{minresthm}.}
Let $\tau':X\to\tilde X$ be the contraction of all the $(-1)$-curves in
$\tau^{-1}(p)$.
We claim that $\tilde X$ is isomorphic to $\bar X$
and $\tau=\tau'\circ\bar\tau$,
namely there is no $(-1)$-curve in $\tau'(\tau^{-1}(p))$.
The only way for a $(-1)$-curve
to be created after blowing down a $(-1)$-curve $F_j=\pi^{-1}(E_j)$ is for a smooth
rational curve
$F_k$ of self-intersection $-2$ to meet the given $(-1)$-curve $F_j$.
Since no two branched curve meet on $Y$ and $E_j$ is branched
by Lemma \ref{Lem10.2},
then $F_k$ must lie over an unbranched curve $E_k$
which meets $E_j$ at  one point.
Therefore $\pi^*(E_k)$ cannot split (its divisor $\Gamma_k$ is not even)
so that $F_k =\pi^*(E_k)$ and $F_k^2 = 2E_k^2$.
So $E_k$ is a $(-1)$-curve on $Y$
and, since $F_k$ is smooth and rational, the divisor $\Gamma_k$
must consist of two simple points (one of which is the intersection
point with $E_j$).
In this case we see that $E_k$ and $E_j$ should be blown down
on $Y$, so that both of these curves were unnecessarily
blown up.
This proves our claim.
Let $F_j=\pi^{-1}(E_j)$ be a $(-1)$-curve of $X$.
Let $E_k$ be an exceptional curve that meets $E_j$.
Since $E_j$ is branched, $E_k$ is unbranched and their intersection
is transversal, as any intersection of the $E_i$'s.
Hence $F_j$ contracts to a smooth point of $F_k$
(cf.\ Lemma \ref{doublecovers}).
Finally, $E_j^2=-2$ implies that $E_j$ meets at most three
exceptional curves of $\sigma$, namely one curve that
corresponds to a blown up point on $E_j$ and
the exceptional curves on which $q_j$ lies, that are at most two.
\qed

\medskip
The previous analysis offers also an alternative route to obtaining
the minimal resolution of $p\in X_0$:
we may first contract the branched $(-2)$-curves among
the $E_i$'s and then take the double cover,
namely the diagram
\[
\begin{array}{ccc}
X & \mapright{\tau'} & \bar X\\
\mapdown\pi & & \mapdown{\bar\pi}\\
Y & \mapright{\sigma'} & \bar Y
\end{array}
\]
commutes, where $\sigma'$ is the contraction of the
branched $(-2)$-curves in $\sigma^{-1}(q_1)\subset Y$
and $X$ is the fiber product $\bar X\times_{\bar Y}Y$.

Note that contracting a $(-2)$-curve on a smooth surface produces
a singularity, namely an ordinary double point (of type $A_1$, see p.\ \pageref{rdp}).
Since the $(-2)$-curves are branched for $\pi:X\to Y$,
the singular points of $\bar Y$ must be considered as
\emph{branched points} for $\bar\pi:\bar X\to\bar Y$.

Let us denote by $\bar F_i$\newcit{defbarF_i}
(resp.\ $\bar E_i$) the image of $F_i$
in $\bar X$ (resp.\ of $E_i$ in $\bar Y$).
For simplicity, suppose that we blow down only a $(-1)$-curve $F_k$
(the general case can be computed inductively).
If $E_k$ meets two unbranched divisors $E_i$ and $E_j$ in $X$,
then $\bar E_i$ and $\bar E_j$ meets in a branched point, hence:
\begin{equation}
\label{barij} \bar F_i\cdot\bar F_j=
\begin{cases}
1 & \text{if $E_i\cdot E_k=E_j\cdot E_k=1$,} \\
F_i\cdot F_j & \text{if $E_i\cdot E_k=0$ or $E_j\cdot E_k=0$.}\\
\end{cases}
\end{equation}

Moreover the self-intersection numbers change as follows:
\begin{equation}\label{barisquare}
\bar F_i^2=\begin{cases} F_i^2+1 & \text {if $E_i\cdot E_k=1$,} \\
F_i^2 & \text{otherwise,}\\
\end{cases}
\end{equation}
while the arithmetic genera stay unchanged.
Formula \fref{fc} and the fact that
$\tau^{{\prime}{-1}}(\bar\tau^{-1}(m_{p,X_0}))=\tau^{-1}(m_{p,X_0})$
imply that the fiber cycle of the minimal resolution
$\bar\tau:\bar X\to X$ is:
\begin{equation}
\label{minfybrecycle}
\bar F =\sum m_{1i}(1+\veps_i)\bar F_i.
\end{equation}

\section{The fundamental cycle}\label{fundcycle}

The fundamental cycle of the (canonical) resolution
is the unique smallest positive cycle:
\begin{equation}
\label{F}
Z=\sum _{i=1}^nz_iF_i
\end{equation}
with $z_i>0$ such that  $Z\cdot F_k\leq0$ for every $k=1,\ldots,n$.
If $F_i$ splits in $F'_i$ and $F''_i$ as in Remark \ref{remFi},
\emph{a priori} we should consider in \fref{F} two distinct
coefficients $z'_i$ and $z''_i$.
But if $z'_i$ were different from $z''_i$ in \fref{F},
then we could exchange them and take the g.c.d.\ of the cycles,
so $z_i=\min\{z'_i,z''_i\}$ would fulfill the properties
of the fundamental cycle.
Therefore we may and will assume $z'_i=z''_i=z_i$.

In general the fundamental cycle of any resolution can be
computed inductively as follows.
Let $F_1,\ldots,F_n$ be the exceptional curves.
\begin{enumerate}
\setlength{\parskip}{0cm}
\setlength{\itemsep}{0cm}

\item[(1)] Set $Z=\sum_{i=1}^n F_i$.

\item[(2)] Check if $Z\cdot F_j\le 0$ for every $j$.

\item[(3)] If (2) is false, there exists $j$ such that $Z\cdot F_j>0$.
Replace $Z$ with $Z+F_j$ and go back to (2).

\item[(4)] Otherwise, if (2) is true, $Z$ is the fundamental cycle.

\end{enumerate}

In \cite{laufer1} Laufer used essentially the properties of the
fundamental cycle in order to describe precisely the relation
between the topological types of the canonical and the minimal
resolution.
However Laufer showed only an implicit formula for $Z$
(see Lemma 5.5 in \cite{laufer1}).
In the next theorem we give an explicit formula for $Z$,
that turns out to be very simple and that may help to
understand better the algorithms described by Laufer
in \cite[Theorems 5.7 and 5.10]{laufer1}.
It is natural to compare the fiber cycle with the fundamental one.
The definitions implies that $F\ge Z$.
It is known that the equality holds for every resolution
for special types of singularities, for example rational \cite{a2}
and minimally elliptic \cite[Theorem 3.13]{laufer2} ones.
Regarding double points,
Dixon showed that $Z=F$ for every resolution if
$\mu_1=\mult_{q_1}(B_0)$ is even
\cite[Theorem 1]{dixon}, and that $Z=F$ for the minimal resolution
if $B_0$ is analytically irreducible at $q_1$ \cite[Theorem 2]{dixon}.
In section \ref{flowers} we will classify all the double point
singularities for which $F>Z$.

\begin{theorem}
\label{numcycle}
The fundamental cycle $Z$ of the canonical resolution
$\tau:X\to X_0$ of $p\in X_0$ differs from the fiber cycle
$F=\pi^*(E_1^*)$ if and only if there exists $j>1$
such that $\veps_j=0$, $q_j$ is infinitely near of order one
to $q_1$ and:
\begin{equation}
\label{A}
m_{1i}+m_{ji} \hbox{ is even for every $i$ such that $\veps_i=0$.}
\end{equation}
In that case the fundamental cycle is:
\begin{equation}
\label{Fcycle}
Z=\frac{1}{2}\pi^*(E_1^*+E_j^*)
=\sum_{i=1}^n\frac{1}{2}(1+\veps_i)(m_{1i}+m_{ji})F_i.
\end{equation}
\end{theorem}

We remark that condition \fref{A} implies that $\veps_1=1$, because
$m_{11}=1$ and $m_{j1}=0$, thus the multiplicity
$\tilde\alpha_1=\mu_1$ of $B_0$ at $q_1$ is odd. Furthermore
$j$ is uniquely determined. Indeed, if $q_i>^1q_1$ (and $i\ne j$),
then $m_{1i}=1$ and $m_{ji}=0$, so \fref{A} would imply that
$\veps_i=1$.

\medskip
\begin{pf}
By \fref{F*.Fi}, $Z\le F$.
If $Z\neq F$, we may write:
\begin{equation}
\label{fzp}
F=Z+P
\end{equation}
where $P=\sum_i t_i F_i$ is a positive (non-zero) divisor.
It follows from \fref{fcsquare} that:
\[
-2=F^2=Z^2+P^2+2Z\cdot P.
\]
Since $Z>0$ and $P>0$, then $Z^2<0$ and $P^2<0$.
Moreover $Z\cdot P\leq 0$ because $Z$ is the fundamental cycle,
so the only possibility is:
\[
Z^2=P^2=-1,\qquad Z\cdot P=0.
\]
Hence $F\cdot P=-1$ and by formulas \fref{F*.Fi}, one finds that:
$
-1=F\cdot P=F\cdot t_1F_1=-(2-\veps_1)t_1,
$
which forces
$
t_1=\veps_1=1,
$
so $E_1$ must be branched and $F_1$ is forced to belong to $P$ with multiplicity one.
Since $Z\cdot F_k\le0$ for each $k$, $Z^2=-1$ implies that there
exists an unique $j$ such that:
\[
Z\cdot F_j=-1,\qquad Z\cdot F_k=0\quad\hbox{for every $k\neq j$,}
\]
and $z_j=1$.
Therefore $t_j=-Z\cdot P=0$, i.e., $F_j$ is not a component of $P$,
and the coefficient of $F_j$ in $F$ is $(1+\veps_j)m_{1j}=z_j+t_j=1$,
so
\[
\veps_j=0 \quad\text{and}\quad m_{1j}=1,
\]
in particular $j\neq 1$.
It follows from \fref{fzp} that:
$
P\cdot F_1=-1
$,
\, $P\cdot F_j=1$,\, $P\cdot F_k=0$\,for $k\neq1,j$.
The previous three equations are equivalent to:
\[
\sum_{i=1}^n (2-\veps_i)t_i s_{ik}=
\begin{cases} -2 & \text{if $k=1$}\\ 1 & \text{if
$k=j$}\\ 0 & \text{if $k\neq1,j$.}
\end{cases}
\]
Recalling that:
\[
\sum_{i=1}^n 2m_{1i} s_{ik}=2(E_1^*\cdot E_k)=
\begin{cases} -2 & \text{if $k=1$}\\ 0 &
\text{if $k\neq 1$}
\end{cases}
\]
and setting $m_1$ and $t$ row vectors with the obvious entries,
one finds that:
\[
((2-\veps)t-2m_1)S=e_k
\]
where $e_k$ is the row vector with the $k$-th entry equal to 1 and $0$
everywhere else.
Multiplying both sides with $S^{-1}$, the vector $((2-\veps)t-2m_1)$ is
the $k$-th row of the matrix $S^{-1}$.
In particular:
$
2t_j-2m_{1j}=0-2=-2
$
is the $(j,j)$-entry in $S^{-1}=-M^tM$.
Therefore:
\[
-2=-\sum_{1\leq i\leq j}m_{ij}^2=-m_{1j}^2-m_{jj}^2-\sum_{1<i<j}
m_{ij}^2=-2-\sum_{1<i<j}m_{ij}^2
\]
that is possible if and only if $m_{ij}=0$ for every $1<i<j$.
This means that $q_j$ is proximate to $q_1$, but $q_j\not>q_i$
for $i\ne1$, i.e.\ $q_j>^1q_1$.
Furthermore, the $(k,j)$-entry in $S^{-1}$ is:
\[
(2-\veps_k)t_k-2m_{1k}=-\sum_{1\leq i\leq j}m_{ik}m_{ij}=-m_{1k}-m_{jk},
\]
that we may rewrite as follows:
$
t_k=\frac{1}{2}(1+\veps_k)(m_{1k}-m_{jk})
$.
Since the coefficient of $F_k$ in $F$ is $(1+\veps_k)m_{1k}=t_k+z_k$, then
\begin{equation}
\label{zk}
z_k=\frac{1}{2}(1+\veps_k)(m_{1k}+m_{jk})
\end{equation}
must be an integer, that proves \fref{A} and \fref{Fcycle}.

Vice versa, if there exists $j$ as in the statement, we may order the
blowing-ups $\sigma_i$ in such a way that $j=2$.
It suffices to show that $Z'=\pi^*(E_1^*+E_2^*)/2$ has the property
that $Z'\cdot F_k\leq 0$ for every $k=1,\ldots,n$.
Indeed $Z'<F$, so the first part of the proof implies that $Z'$ has
to be the fundamental cycle.
If $k>2$, then:
$
Z'\cdot F_k=(1-\veps_k/2)(E_1^*\cdot E_k+E_2^*\cdot E_k)=0.
$
Moreover $Z'\cdot F_2=E_2^{*2}=-1$ and
$Z'\cdot F_1=E_1^*\cdot E_1/2+E_2^*\cdot E_1/2=0$.
\qed
\end{pf}

If $F>Z$, the arithmetic genus of $Z$ is, by $Z^2=-1$,
\fref{Fcycle} and \fref{as}:
\begin{eqnarray}
p_a(Z)&=&\frac{E_1^*+E_j^*}{2}\cdot \left(\sigma^*\left(K_{Y_0}+\frac{B_0}{2}\right)
-\sum_{k=1}^n \left(\frac{\alpha_k}{2}-1\right)E_k^*\right)+\frac{1}{2}=\nonumber\\
&=&\frac{\alpha_1+\alpha_j-2}{4}=\frac{\tilde\alpha_1+\tilde\alpha_j-2}{4}.
\label{numcyclegenus}
\end{eqnarray}
where the last equality follows from $\alpha_1=\tilde\alpha_1-1$ and
$\alpha_j=\tilde\alpha_j+1$.

\medskip
Now we compute the fundamental cycle of the minimal resolution.

\begin{lemma}
Let $\bar\tau:\bar X\to X_0$ be the minimal resolution of $p\in X_0$.
Then the fundamental cycle $\bar Z$ of $\bar\tau$ is:
\begin{equation}\label{barF}
\bar Z=\sum_iz_i\bar F_i
\end{equation}
where $\bar F_i$ are the exceptional curves of $\bar\tau$ and $Z=\sum_i z_i F_i$ is
the fundamental cycle of the canonical resolution.
\end{lemma}

\begin{pf}
Without any loss of generality, we may assume to blow down only
a $(-1)$-curve $F_k=\pi^{-1}(E_k)$, where $E_k^2=-2$ and $\veps_k=1$.
Recall that $E_k$ meets at least one and at most three unbranched divisors.
First, we shall prove that:
\begin{equation}\label{barFi.Fj}
\sum_{i\ne k}z_i\bar F_i\cdot \bar F_j\le0.
\end{equation}
for every $j\ne k$.
For this purpose, we claim that:
\begin{equation}\label{rk}
z_k=\sum_{i:E_i\cdot E_k=1}z_i.
\end{equation}
Suppose that $E_k$ meets three unbranched divisors, i.e.,
$q_k=E_{k_1}\cap E_{k_2}$
and we blew up a point $q_{k_3}$ lying on $E_k$,
with $\veps_{k_i}=0$ for $i=1,2,3$.
Then $m_{1k}=m_{1k_3}=m_{1k_1}+m_{1k_2}$, so:
\[
2m_{1k}=(1+\veps_k)m_{1k}=\sum_{i=1}^3(1+\veps_{k_i})m_{1k_i}
=m_{1k_1}+m_{1k_2}+m_{1k_3}
\]
which proves \fref{rk} if $Z=F$, by \fref{fc}.
Similarly, $2m_{2k}=m_{2k_1}+m_{2k_2}+m_{2k_3}$ and \fref{rk} holds
even if $Z<F$, by \fref{zk}.
If $E_k$ meets only one or two unbranched divisors, there are four
possible configurations and the proof of \fref{rk} is analogous.

Clearly \fref{barFi.Fj} holds if $E_j\cdot E_k=0$.
Otherwise if $E_j\cdot E_k=1$, by \fref{barij},
\fref{barisquare} and the fact that if $E_j\cdot E_k=E_i\cdot E_k=1$
then $F_i\cdot F_j=0$, formula \fref{rk} implies that:
\begin{align*}
\sum_{i\ne k}z_i\bar F_i\cdot\bar F_j &= \sum_{i:E_i\cdot
E_k=1}z_i\bar F_i\cdot\bar F_j+
\sum_{i:E_i\cdot E_k=0}z_iF_i\cdot F_j=\\
&=\sum_{i:E_i\cdot E_k=1}z_i+\sum_{i\ne k}z_iF_i\cdot F_j
=z_k+\sum_{i\ne k}z_iF_i\cdot F_j=Z\cdot F_j\le0
\end{align*}
which proves \fref{barFi.Fj}.
Let $\bar Z=\sum_{i\ne k}s_i\bar F_i$ be the fundamental cycle
of $\bar\tau$.
If we  show that for every $j$:
\begin{equation}\label{F+Fk}
\left(\sum_{i\ne k}s_iF_i+\sum_{i:E_i\cdot E_k=1}s_iF_k\right)\cdot F_j\le0,
\end{equation}
then $s_i=z_i$, for $i\ne k$, and \fref{barF} holds. Indeed if
$E_j\cdot E_k=0$ then \fref{F+Fk} is trivial. If $E_j\cdot E_k=1$,
then the left hand side of \fref{F+Fk} becomes:
\[
\sum_{i:E_i\cdot E_k=0}s_iF_i\cdot F_j+\sum_{i:E_i\cdot E_k=1}s_iF_i\cdot F_j
+\sum_{i:E_i\cdot E_k=1}s_iF_k\cdot F_j
=\sum_{i\ne k}s_i\bar F_i\cdot\bar F_j=\bar F\cdot\bar F_j\le0.
\]
Finally, for $j=k$, the left hand side of \fref{F+Fk} is:
$$
\sum_{i\ne k}s_iF_i\cdot F_k+\sum_{i:E_i\cdot E_k=1}s_iF_k^2
=\sum_{i:E_i\cdot E_k=1}s_iF_i\cdot F_k
-\sum_{i:E_i\cdot E_k=1}s_i=0.\eqno\Box
$$
\end{pf}

\begin{corollary}\label{corfundcycle}
Let $F,Z$ (resp.\ $\bar F,\bar Z$) be the fiber and the
fundamental cycle of the canonical (resp.\ minimal) resolution.
Then $F>Z$ and $\bar F=\bar Z$ if and only if $q_2$ is the unique
proximate point to $q_1$ and $\tilde\alpha_1=\tilde\alpha_2$ is odd.
Furthermore, this happens if and only if $\bar F^2=\bar Z^2=-1$.
\end{corollary}

\begin{pf}
Suppose that $F>Z$ and $\bar F=\bar Z$.
This means that $F_1$ is a $(-1)$-curve that we blow down,
hence $E_1^2=-2$ and there is only one proximate point to $q_1$,
that is $q_2$, so $\tilde\alpha_1=\tilde\alpha_2$.
Moreover $\tilde\alpha_1$ is odd by Theorem \ref{numcycle}.
Conversely, if $E_1$ is branched and $E_1^2=-2$,
then $F_1$ is a $(-1)$-curve that we blow down.
Hence $m_{2i}=m_{1i}$ for every $i>2$, therefore the coefficient
of $F_i$ in $F$ is the same as the
coefficient of $F_i$ in $Z$, for every $i>2$.
The last assertion follows from the fact that if
the fundamental cycle of a resolution of $p$ has self-intersection
$-2$, then on any resolution the fundamental cycle is equal to
the fiber cycle (cf.\ \cite[Lemma 5.2]{laufer3}
or \cite[p.\ 110]{dixon}).
\qed
\end{pf}

\section{The description of the Enriques digraph}\label{flowers}

We want to describe the weighted Enriques digraph
of those double point singularity for which
the fundamental cycle of the canonical resolution
is strictly contained in the fiber cycle.
Recall that the weight of the vertex $q_i$ is $\mu_i$,
i.e.\ the multiplicity at $q_i$
of the branch curve $B_{i-1}$ of $\pi_{i-1}:X_{i-1}\to Y_{i-1}$.
Before going on, we need some remark about proximate points.

Let us call \emph{proximity subgraph} of $q_1$
the subgraph of the Enriques digraph consisting only
of the proximate points to $q_1$ (and the arrows among them).

\begin{figure}[ht]
\centerline{
\begin{picture}(220,180)
\thinlines \put(130,77){\line(0,-1){54}}
\put(130,83){\line(0,1){54}} \put(127,80){\line(-1,0){54}}
\put(128.66,77.32){\line(-1,-2){27.32}}
\put(127.9,77.9){\line(-1,-1){55.8}}
\put(127.32,78.66){\line(-2,-1){114.64}}
\put(127.32,81.34){\line(-2,1){54.64}}
\put(128.336,82.496){\line(-2,3){56.672}}
\put(131.34,82.68){\line(1,2){27.32}}
\put(132.496,81.664){\line(3,2){85.008}}
\put(103,20){\line(1,0){24}} \put(73,20){\line(1,0){24}}
\put(55,20){\line(1,0){12}} \put(13,20){\line(1,0){10}}
\put(70,83){\line(0,1){24}} \put(70,113){\line(0,1){10}}
\put(70,155){\line(0,1){12}} \put(133,140){\line(1,0){24}}
\put(163,140){\line(1,0){10}} \put(205,140){\line(1,0){12}}
\put(10,20){\circle{6}} \put(70,20){\circle{6}}
\put(100,20){\circle{6}} \put(130,20){\circle{6}}
\put(130,80){\circle{6}} \put(70,80){\circle{6}}
\put(70,110){\circle{6}} \put(70,170){\circle{6}}
\put(130,140){\circle{6}} \put(160,140){\circle{6}}
\put(220,140){\circle{6}} \put(35,20){\circle*{2}}
\put(40,20){\circle*{2}} \put(45,20){\circle*{2}}
\put(70,135){\circle*{2}} \put(70,140){\circle*{2}}
\put(70,145){\circle*{2}} \put(185,140){\circle*{2}}
\put(190,140){\circle*{2}} \put(195,140){\circle*{2}}
\put(134,75){$\scriptstyle q_1$} \put(127,10){$\scriptstyle q_2$}
\put(97,10){$\scriptstyle q_3$} \put(67,10){$\scriptstyle q_4$}
\put(6,10){$\scriptstyle q_{i_1}$} \put(45,79){$\scriptstyle
q_{i_1+1}$} \put(45,109){$\scriptstyle q_{i_1+2}$}
\put(55,169){$\scriptstyle q_{i_2}$} \put(120,148){$\scriptstyle
q_{i_2+1}$} \put(150,148){$\scriptstyle q_{i_2+2}$}
\put(217,148){$\scriptstyle q_{i_3}$} \put(160,60){\small(etc.)}
\thicklines \put(118,20){\vector(1,0){0}}
\put(88,20){\vector(1,0){0}} \put(58,20){\vector(1,0){0}}
\put(28,20){\vector(1,0){0}} \put(70,92){\vector(0,-1){0}}
\put(70,122){\vector(0,-1){0}} \put(70,152){\vector(0,-1){0}}
\put(142,140){\vector(-1,0){0}} \put(172,140){\vector(-1,0){0}}
\put(202,140){\vector(-1,0){0}} \put(130,53){\vector(0,1){0}}
\put(117.925,55.85){\vector(1,2){0}}
\put(110.908,60.908){\vector(1,1){0}}
\put(105.85,67.925){\vector(2,1){0}} \put(103,80){\vector(1,0){0}}
\put(105.85,92.075){\vector(2,-1){0}}
\put(115.023,102.4654){\vector(2,-3){0}}
\put(130,107){\vector(0,-1){0}}
\put(142.075,104.15){\vector(-1,-2){0}}
\put(152.4654,94.977){\vector(-3,-2){0}}
\end{picture}}
\caption{The proximity subgraph of $q_1$}
\label{figuregraph}
\end{figure}
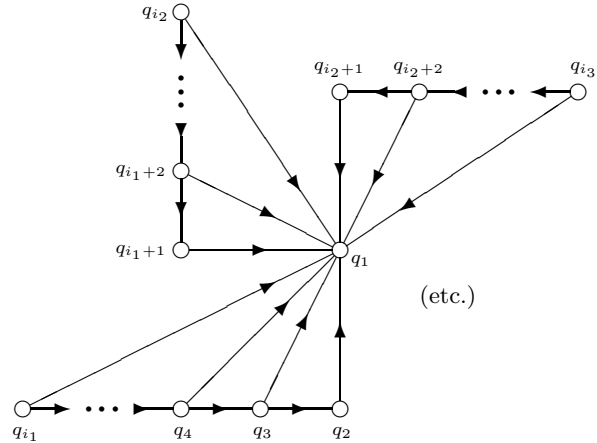

We may order the $\sigma_i$'s (the blowing-ups) in such a way
that $q_2$, $\ldots$, $q_{n'}$ are all the proximate points to $q_1$
and for every $j=1,\ldots,h$:
\begin{equation}\label{qij}
q_{i_j}>^1q_{i_j-1}>^1\cdots>^1q_{i_{j-1}+2}>^1q_{i_{j-1}+1}>^1q_1,
\end{equation}
where $n'=i_h>i_{h-1}>\cdots>i_2>i_1>i_0=1$.
Thus the proximity digraph of $q_1$ has the shape of Figure
\ref{figuregraph}, which looks like a \emph{flower} with $h$ petals.
We say that the $j$-th petal has \emph{length} $i_j-i_{j-1}$.

Clearly the proximity subgraph of any point has a similar shape.

Let us say that a vertex of the Enriques digraph is
\emph{very odd} if its weight is odd, its proximity subgraph has exactly one petal
of odd length and all the other petals of even length.
Now we are ready to prove the following:

\begin{theorem}\label{thmgraph}
The fundamental cycle is strictly contained in the fiber cycle of
the canonical resolution of $p\in X_0$ if and only if
the weighted Enriques digraph of $q_1=\pi(p)$
has the following properties:
\begin{enumerate}
\setlength{\parskip}{0cm}
\setlength{\itemsep}{0cm}

\item $q_1$ is a very odd vertex (in particular its weight $\mu_1$ is odd);

\item a proximate point $q_i$ to $q_1$,
belonging to a petal of even
(resp.\ odd) length of the proximity subgraph of $q_1$,
is a very odd vertex if and only if $q_i$ is
infinitely near of odd (resp.\ even) order to $q_1$;

\item inductively, 1 and 2 hold replacing $q_1$ with any very odd vertex.
\end{enumerate}
\end{theorem}

\begin{pf}
Suppose that the fundamental cycle $Z$ is strictly contained
in the fiber cycle $F$.
By Theorem \ref{numcycle}, $\mu_1$ is odd and there exists $j$
such that $\veps_j=0$, $q_j>^1q_1$ and condition \fref{A} holds.
Moreover we may and will assume that $j=2$,
so $\veps_2=0$ and $\mu_2$ is even.

Consider the proximity subgraph of $q_1$ as above
(cf.\ formula \fref{qij} and Figure \ref{figuregraph}).
We claim that $i_j$ is even for every $j=1,\ldots,h$,
thus $q_1$ is a very odd vertex.

Suppose that $i_1>2$, namely in the canonical resolution process
we blow up $q_3=E_2\cap E_1$. Then $m_{13}=2$ and $m_{23}=1$,
so condition \fref{A} implies that $\veps_3=1$ and $\mu_3$ is odd.
Since $\veps_1=\veps_3=1$, the intersection $q_4=E_3\cap E_1$ is a
singular point of the branch curve $B_3$ of $\pi_3:X_3\to Y_3$, so
we must blow up also $q_4$. Similarly, if we blow up $q_5=E_4\cap
E_1$, then $m_{15}+m_{25}=5$, thus condition \fref{A} forces
$\veps_5=1$ and we must blow up also $q_6=E_5\cap E_1$. Repeating
this argument, it follows that the first petal has odd length and
$i_1$ is even.

Look at the second petal. Now $q_{i_1+1}>^1q_1$, so
$m_{1,i_1+1}=1$ and $m_{2,i_1+1}=0$. Hence \fref{A} implies that
$\veps_{i_1+1}=1$ and $\mu_{i_1+1}$ is odd. Therefore we must blow
up $q_{i_1+2}=E_{i_1+1}\cap E_1$. If $i_2>i_1+2$,
it means that we blow up also $q_{i_1+3}=E_{i_1+2}\cap
E_1$, then $m_{1,i_1+3}+m_{2,i_1+3}=3$, so \fref{A} forces
that $\veps_{i_1+3}=1$ and we must blow up $q_{i_1+4}=E_{i_1+3}\cap
E_1$ too. Proceeding in this way, this shows that the second
petal has even length and $i_2$ is even.
The same argument works for the $j$-th petal, with $j>2$,
just by replacing $i_1$ with $i_{j-1}$.
This proves our claim that $q_1$ is a very odd vertex.

Now we want to show that $\veps_i\equiv i\pmod2$
for every $i=1,\ldots,i_h$, and
$q_{2l-1}$ is a very odd vertex for every $l=2,\ldots,i_h/2$.

We already know that $\veps_2=0$ and $\veps_{2i-1}=1$ for every
$i=1,\ldots,i_h/2$. Suppose by contradiction that $\veps_4=1$ (and
$i_h>2$). Since $\veps_3=1$, we must blow up also $q_k=E_4\cap E_3$
and $m_{1k}+m_{2k}$ is odd, so condition \fref{A} implies that
$\veps_k=1$. Hence we must blow up $q_{k+1}=E_k\cap E_4$ too, and
$m_{1,k+1}+m_{2,k+1}$ is again odd, thus $\veps_{k+1}=1$ by \fref{A}.
Going on in this way, we produce each time another branched
exceptional curve, so we should never stop blowing up,
contradicting Theorem \ref{thmcanres}.
This shows that $\veps_4=0$ and $\mu_4$ is even.
The proof that $\veps_{2l}=0$ for every $l=3,\ldots,i_h/2$
is similar.

Consider the proximity subgraph of $q_{2l-1}$, for
$l=2,\ldots,i_h/2$. Repeating exactly the same arguments as
for $q_1$ and $q_2$,
one finds out that $q_{2l}$ (which is proximate to $q_{2l-1}$)
belongs to a petal of odd length, while all other petals of the
proximity subgraph of $q_{2l-1}$ have even length,
thus $q_{2l-1}$ is a very odd vertex.

It remains to prove that the proximity subgraph of a non-very-odd vertex
can be arbitrary.
For every $l=1,\ldots,i_h/2$, we proved that $\veps_{2l}=0$,
so $q_{2l}$ cannot be very odd.
Moreover $m_{1,2l}+m_{2,2l}$ is even.
If $q_k$ is proximate to $q_{2l}$ (and $k\ne 2l+1$),
then $m_{1k}+m_{2k}$ is a multiple of $m_{1,2l}+m_{2,2l}$, hence it
is even and \fref{A} imposes no condition on $q_k$.
This means that the proximity subgraph of a non-very-odd vertex, as $q_{2l}$,
can be arbitrary and concludes the proof
that the Enriques digraph has properties 1, 2 and 3.

Conversely, suppose that the three properties hold.
One may easily check that the $m_{ij}$'s satisfy
condition \fref{A}, where the wanted $q_j$ is the infinitely
near point of order one to $q_1$ belonging to the petal of
odd length, therefore one concludes by Theorem \ref{numcycle}.
\qed
\end{pf}

Note that, with the notation of the proof,
$\tilde\alpha_1$ and $\tilde\alpha_2$ are odd, while
$\tilde\alpha_i$ is even for every $i=3,\ldots,i_h$, by \fref{alfa'}.
Moreover $\veps_1=1$ forces $\tilde B\cdot E_1=0$, or equivalently
$\tilde\alpha_1=\sum_{j=1}^n\tilde\alpha_jq_{1j}=\sum_{j=2}^{i_h}\tilde\alpha_j$.
By induction on the number $i_h$ of proximate points to $q_1$,
it is easy to check that $\smash{\tilde\alpha_1=\tilde\alpha_2+\sum_{j=3}^{i_h}\tilde\alpha_j
\equiv\tilde\alpha_2\pmod4}$,
thus $\tilde\alpha_1+\tilde\alpha_2\equiv2\pmod4$
(cf.\ the genus formula \fref{numcyclegenus}).

\section{Some examples}

\begin{example}\label{ex1}
Let $B_0$ be defined by: $y(y-x^2)(y+x^2)=0$. One usually says
that $B_0$ has two infinitely near triple points at $q_1$. Our
combinatorial data are:
\[
\raisebox{-7pt}{\begin{picture}(60,20)
\thinlines
\put(14.5,10){\line(1,0){31}}
\put(10,10){\circle{9}}
\put(50,10){\circle{9}}
\put(8,7.8){$\scriptstyle3$}
\put(48,7.8){$\scriptstyle4$}
\thicklines
\put(27,10){\vector(-1,0){0}}
\end{picture}}
\qquad\hbox{or equivalently }\qquad M=
\begin{pmatrix}
1 & 1 \\
0 & 1
\end{pmatrix},
\quad
\begin{array}{l}
\mu=(3,4), \\
\veps=(1,0).
\end{array}
\]
\end{example}

The exceptional curves for $\tau:X\to X_0$ are
a smooth rational curve $F_1$ with $F_1^2=-1$
and a smooth elliptic curve $F_2$ with $F_2^2=-2$
that meet in a point $P$.
By Theorems \ref{fibercycle} and \ref{numcycle},
the fiber cycle of the canonical resolution is $F=2F_1+F_2$,
while the fundamental cycle is $Z=F_1+F_2<F$,
as one may also check directly.
Moreover ${F_2}_{|F_2}=-2P$ by \fref{Fi|Fi}.

The minimal resolution $\bar\tau:\bar X\to X_0$ is obtained by
contracting the $(-1)$-curve $F_1$.
Therefore $\bar F_2$ is the only exceptional curve for $\bar\tau$
and $\bar F_2$ is a smooth elliptic curve with $\bar F_2^2=-1$.
Clearly the fiber cycle and the fundamental cycle of the
minimal resolution $\bar\tau$ are $\bar Z=\bar F=\bar F_2$.

The previous example can be generalized as follows.

\begin{example}\label{ex2}
Let $B_0$ be a curve with $2k$ infinitely near points
$q_1$, $\ldots$, $q_{2k}$ of the same odd multiplicity
$\tilde\alpha_1=\cdots=\tilde\alpha_{2k}=2g+1$, for some $g\ge1$.
More precisely, $q_i>^1q_{i-1}$ for $1<i\le2k$
and the weighted Enriques digraph is:
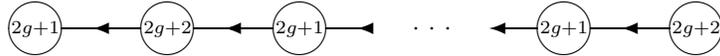
\[
\begin{picture}(280,20)
\thinlines
\put(20,10){\line(1,0){30}}
\put(70,10){\line(1,0){30}}
\put(120,10){\line(1,0){15}}
\put(185,10){\line(1,0){15}}
\put(220,10){\line(1,0){30}}
\put(10,10){\circle{19}}
\put(60,10){\circle{19}}
\put(110,10){\circle{19}}
\put(210,10){\circle{19}}
\put(260,10){\circle{19}}
\put(154,10){\circle*{1}}
\put(160,10){\circle*{1}}
\put(166,10){\circle*{1}}
\put(1,7.8){$\scriptstyle 2g+1$}
\put(51,7.8){$\scriptstyle 2g+2$}
\put(101,7.8){$\scriptstyle 2g+1$}
\put(201,7.8){$\scriptstyle 2g+1$}
\put(251,7.8){$\scriptstyle 2g+2$}
\thicklines
\put(32,10){\vector(-1,0){0}}
\put(82,10){\vector(-1,0){0}}
\put(132,10){\vector(-1,0){0}}
\put(182,10){\vector(-1,0){0}}
\put(232,10){\vector(-1,0){0}}
\end{picture}
\]
\end{example}

The exceptional curves for $\tau:X\to X_0$ are the following:
$(-1)$-curves $F_{2i-1}$, for every $i=1,\ldots,k$;
smooth rational curves $F_{2i}$ with
self-intersection $-4$, for $i=1,\ldots,k-1$,
and a smooth curve $F_{2k}$ of genus $g$ with $F_{2k}^2=-2$.
By Theorems \ref{fibercycle} and \ref{numcycle},
the fiber cycle and the fundamental
cycle of the canonical resolution are respectively:
\[
F=\sum_{i=1}^k(2F_{2i-1}+F_{2i}),
\qquad Z=F_1+F_2+\sum_{i=2}^k(2F_{2i-1}+F_{2i}).
\]
Blow down the $F_{2i-1}$'s, for $i=1,\ldots,k$, thus
the exceptional curves for the minimal resolution
$\bar\tau:\bar X\to X_0$ are the $\bar F_{2i}$'s, for $i=1,\ldots,k$,
which are smooth rational curves with self-intersection $-2$,
except $\bar F_{2k}$ which is smooth of genus $g$ with $\bar F_{2k}^2=-1$.
The fundamental cycle equals the fiber cycle of the minimal resolution
$\bar Z=\bar F=\sum_{i=1}^k\bar F_{2i}$.

\begin{example}(cf.\ \cite[p.~322]{laufer3})\label{ex3}
Let $B_0$ be defined by: $y(x^4+y^6)=0$.
In this case our combinatorial data are:
\[
\raisebox{-25pt}{\begin{picture}(60,60)
\thinlines
\put(14.5,10){\line(1,0){31}}
\put(10,14.5){\line(0,1){31}}
\put(14.5,50){\line(1,0){31}}
\put(13.5,13.5){\line(1,1){33}}
\put(10,50){\circle{9}}
\put(50,10){\circle{9}}
\put(10,10){\circle{9}}
\put(50,50){\circle{9}}
\put(48,7.8){$\scriptstyle2$}
\put(8,7.8){$\scriptstyle5$}
\put(8,47.8){$\scriptstyle3$}
\put(48,47.8){$\scriptstyle4$}
\thicklines
\put(27,10){\vector(-1,0){0}}
\put(10,27){\vector(0,-1){0}}
\put(27,50){\vector(-1,0){0}}
\put(27,27){\vector(-1,-1){0}}
\end{picture}}
\qquad
M=\begin{pmatrix}
1 & 1 & 1 & 2\\
0 & 1 & 0 & 0 \\
0 & 0 & 1 & 1 \\
0 & 0 & 0 & 1
\end{pmatrix},
\quad
\begin{array}{l}
\mu=(5,2,3,4),\\
\veps=(1,0,1,0).
\end{array}
\]
\end{example}

The fiber cycle of the canonical resolution is
$F=2F_1+F_2+2F_3+2F_4$, while the fundamental cycle is
$Z=F_1+F_2+F_3+F_4<F$, as it should be by Theorem \ref{thmgraph},
because $q_1$ is a very odd vertex, its proximity digraph has just two
petals (one of length 1 and the other of length 2) and $q_3$ is
also very odd. The minimal resolution is obtained by blowing
down $F_3$, therefore:
\[
\bar F=2\bar F_1+\bar F_2+2\bar F_4>\bar F_1+\bar F_2+\bar F_4=\bar Z.
\]

\medskip\noindent{\it Rational double points}
(see \cite{a2}, \cite{d} and \cite{bpv}).\newcit{rdp}
It is very well-known that the rational double points are given by the
following equations:
\begin{alignat*}{2}
A_n : &\quad z^2=x^2+y^{n+1},  &&\qquad(n\ge1),\\
D_n : &\quad z^2=y(x^2+y^{n-2}), &&\qquad(n\ge3),\\
E_6 : &\quad z^2=x^3+y^4,\\
E_7 : &\quad z^2=x(x^2+y^3),\\
E_8 : &\quad z^2=x^3+y^5
\end{alignat*}
and the minimal resolution consists in smooth rational curves of
self-intersection $-2$ whose dual graph is
the corresponding Dynkin diagram.

Note that, for a rational double point, the fundamental cycle $Z$ of
the canonical resolution equals the fiber cycle $F$.
Indeed if $Z<F$, then formula \fref{numcyclegenus} says that
$p_a(Z)=(\tilde\alpha_1+\tilde\alpha_2-2)/4$, where
$\tilde\alpha_1$ is odd and $\tilde\alpha_1\ge2$, so $p_a(Z)>0$,
contradicting Artin's criterion (which says that $p_a(Z)=0$
if and only if the singularity is rational,
cf.\ \cite[Theorem 3]{a2}).

Moreover,
starting from
the well-known formula for the arithmetic genus of a sum of two curves,
that is
$p_a (C+D)= p_a(C)+ p_a (D) + (C\cdot D) - 1$, and using \fref{isquare-ij}
and \fref{igenus}
we find out that
(cf.\ \fref{fibcyclegenus}):
\[
p_a(Z)=p_a(F)=\frac{1}{2}\sum_{i=1}^n m_{1i}(\gamma_i-(\veps_i-2)E_i^2-4).
\]
Hence we can see directly that $p_a(Z)=0$ if and only if
every branched exceptional divisor has self-intersection $-4$,
every unbranched exceptional divisor has self-intersection $-1$
and $\gamma_i = 2$, or self-intersection $-2$ and $\gamma_i=0$.
Thus every $F_i$ (if $F_i$ splits, every irreducible component of $F_i$)
is rational with self-intersection $-2$ and the
canonical resolution is minimal.

\section{Adjunction conditions}\label{adjoints}

We want to study the conditions that a double
point singularity $p\in X_0$ imposes to canonical and
pluricanonical systems of a surface.
Recall that locally $p$ is $\pi_0^{-1}(q_1)$,
where $\pi_0:X_0\to Y_0$ is a double cover,
$Y_0$ is a smooth surface, $X_0$ is normal and
$q_1$ is an isolated singular point
of the branch curve $B_0$ of $\pi_0$.
Then we consider the canonical resolution $\pi:X\to Y$,
that is a double cover branched along the smooth curve $B$.

In \fref{cd} we defined the adjunction condition divisor $D^*$
as the pullback of
\begin{equation}\label{D}
D=\sum_i (\alpha_i/2-1)E_i^*=\sigma^*(K_{Y_0}+B_0/2)-(K_Y+B/2).
\end{equation}
So it suffices to understand what are the conditions that $D$ imposes
to the adjoint linear system $|K_Y+B/2|$.
It is well known that $D=0$, or equivalently $D^*=0$,
if and only if $p\in X_0$ is a rational double point (cf.\ previous section).

By applying $\sigma_*$ to the exact sequence
$
0 \to \CO_Y(-D) \to \CO_Y \to \CO_{D} \to 0,
$
one sees that $\CI_\Gamma:=\sigma_*\CO_Y(-D)$ is the ideal sheaf
of a zero-dimensional scheme $\Gamma$ supported
at $q_1\in Y_0$.
Let us call $\CI_\Gamma$ the \emph{adjoint ideal} of the singularity.

For our convenience, let us assume that $X_0$ is a double plane,
i.e., $Y_0=\Q$.
Indeed the double point singularity is locally given
as a double cover of an open disc, thus we may always find
an irreducible plane curve $B_0$ of arbitrarily high (and even) degree
whose germ at $q_1$ is analitically isomorphic to
the germ of the branch curve of the double cover at $q_1$.

By \fref{piOx}, \fref{Riemann-H} and the projection formula we have that $\pi_*K_X\cong
K_Y\oplus(K_Y+B/2)$, so $p_g(X)=h^0(X,K_X)=h^0(Y,K_Y+B/2)$ and
$q(X)=h^1(X,K_X)=h^1(Y,K_Y+B/2)$. Riemann--Roch Theorem for $K_Y+B/2$ on
$Y$ and for $K_\Q+B_0/2$ on $\Q$ implies that:
\begin{align*}
&h^0(K_X)-h^1(K_X)=B\cdot(K_Y+B/2)/4+1, \\
&h^0(K_{X_0})=h^0(K_{\Q}+B_0/2)=B_0\cdot(K_\Q+B_0/2)/4+1.
\end{align*}
It follows from \fref{as} that
\begin{equation}
\label{ac}
h^0(K_{X_0})-h^0(K_X)+h^1(K_X)=\sum_{i=1}^n\frac{\alpha_i(\alpha_i-2)}{8}=:c,
\end{equation}
where $c$ is defined by \fref{ac}.
Let us recall a well-known theorem of De Franchis:

\begin{theorem}[De Franchis]
Let $\pi_0:X_0\to\Q$ be a double plane and $\pi:X\to Y$ its canonical
resolution.
Then $q(X)>0$ if and only if there is a plane curve $B'$
(possibly $B'=0$) such that:
\begin{equation}
B_0+2B'=C_1+C_2+\cdots+C_m
\end{equation}
where $C_1,\ldots,C_m$ are curves belonging
to one and the same pencil and $m=2q(X)+2$ (resp.\ possibly $m=2q(X)+1$
if the pencil contains a double curve).
\end{theorem}

\begin{pf}
See \cite{defranchis} (or \cite{cc} for a modern proof).
\qed
\end{pf}

\begin{corollary}\label{adjuncond}
With the above notation, the number $h^0(K_{X_0})-h^0(K_X)$ of
conditions that the singularity $p\in X_0$ imposes to the
canonical system is:
\begin{equation}\label{c=}
c=\sum_{i=1}^n\frac{\alpha_i(\alpha_i-2)}{8}
=\sum_{i=1}^n\frac{\alpha_i/2(\alpha_i/2-1)}{2}.
\end{equation}
\end{corollary}

\begin{pf}
Since we assumed $B_0$ to be irreducible, then $q(X)=h^1(K_X)=0$
by De Franchis' Theorem and \fref{c=} follows from \fref{ac}. \qed
\end{pf}

We remark that De Franchis' Theorem allows us to compute the
adjunction conditions even if $B_0$ were a given reducible curve
and its degree were not assumed to be arbitrarily high.

\medskip
We want to determine which singularity the general element $C$
in $|\CI_\Gamma(h)|=|\sigma_*\CO_Y(K_Y+B/2)|=
\left|\sigma_*\left(\sigma^*(\CO_\Q(hL))\otimes\CO_Y(-D)\right)\right|$
has at $q_1$, where $L$ is
a general line in $\Q$ and
$h=\deg(B_0)/2-3$.
According to formulas \fref{D} and \fref{c=}, one
might expect that $C$ has exactly multiplicity $\alpha_i/2-1$ at $q_i$, for
every $i=1,\ldots,n$. The next example shows that this is not
always the case.

Suppose that $q_1\in B_0$ is the same singularity of Example
\ref{ex1}. Since $\alpha=\mu-\veps$, we have $\alpha_1=2$ and
$\alpha_2=4$, thus one expects the general
element $C$ in $|\CI_\Gamma(h)|$ to pass simply through $q_2$ and
not to pass through $q_1$. But this is not possible, because $q_2$
is infinitely near to $q_1$.

Actually, we see that $(K_Y+B/2)\cdot E_1<0$ and $E_1$ is a fixed
component of $|K_Y+B/2|=|\sigma^*(hL)-E_2|$. Moreover
\begin{equation}\label{K+A-E1}
K_Y+B/2-E_1=\sigma^*(hL)-1\cdot E_1^*-0\cdot E_2^*
\end{equation}
meets non negatively $E_1$ and $E_2$, so $|K_Y+B/2-E_1|$
has no fixed components by the next Lemma \ref{Lultimo}.
Formula \fref{K+A-E1} means that $C$
passes simply through $q_1$ and does not pass through $q_2$ (which
is 1 adjunction condition as well).

\begin{lemma}\label{Lultimo}
Let $\L$ be a linear system on $Y$ which we write as:
\[
\L= \left|\,\sigma^*(hL)-\sum_{i=1}^{n}m_iE_i^* \right|,
\]
where $L$ is a general line in $\Q$, $m_i$ are non-negative
integers and $h$ is arbitrarily high.
Suppose that $\deg\L_{|E_i}\ge0$, for every $i=1,\ldots,n$.
Then $\L$ has no fixed component.
In particular the general member of $\sigma_*\L$ is a plane
curve with multiplicity exactly $m_i$ at $q_i$,
for $i=1,\ldots,n$.
\end{lemma}

\begin{pf}
Since $h\gg 0$, we may assume that the only possible fixed
components of $\L$ are among the $E_i$'s.
For every $i=1$, $\ldots$, $n$,
consider the exact sequence
$
0 \to \L(-E_i) \to \L \to \L_{|E_i} \to 0.
$
We need to show that $h^0(\L(-E_i))<h^0(\L)$.
This will follow from $H^1(\L(-E_i))=0$, because
$H^0(\L_{|E_i})\ne0$ by assumption.
We claim that $R^1\sigma_*\L(-E_i)=0$.
This will imply that $H^1(\L(-E_i))=H^1(\sigma_*\L(-E_i))=0$,
where the last equality follows from Serre's Theorem,
because $h\gg0$, and we will be done.
Indeed $H^1(\L(-E_i)_{|E_1^*})=0$, since
$E_1^*$ is 1-connected, $p_a(E_1^*)=0$
and $\deg\L(-E_i)_{|E_1^*}\ge m_1\ge0$.
\qed
\end{pf}

The above discussion suggested us to introduce the following notion:
we say that a point $q_i$ is \emph{defective}\newcit{defdef}
if there exists a point $q_j$ such that $\alpha_j>\alpha_i$ and $q_j$
is infinitely near of order one to $q_i$.
Hence, if $q_i$ is defective, then $B\cdot E_i<0$,
while we know that $\tilde B\cdot E_i\ge0$ for every $i=0, \ldots, n$,
because $\tilde B$ is the proper transform of a plane curve.
In Example \ref{ex1} (recalled before the previous lemma),
$q_1$ is defective, because $q_2>^1q_1$ and $4=\alpha_2>\alpha_1=2$.

\begin{lemma}\label{D.Ei}\label{D.Ei2}\label{Ldefective}
A point $q_i$ is defective if and only if $D\cdot E_i>0$.
More precisely,
a point $q_i$ is defective if and only if $\veps_i=1$ and there exists
a (necessarily unique) point $q_j>^1q_i$ with
$\tilde\alpha_j=\tilde\alpha_i$ and $\veps_j=0$.
Furthermore, either:
\begin{enumerate}
\setlength{\parskip}{0cm}
\setlength{\itemsep}{0cm}

\item[(i)] $\alpha_i=\tilde\alpha_i-1$, or

\item[(ii)] $\alpha_i=\tilde\alpha_i$ and both $q_i,q_j$ are proximate
to a point $q_k$ with $\veps_k=1$.
\end{enumerate}
Finally a point $q_i$ is defective if and only if $F_i$ is a $(-1)$-curve.
\end{lemma}

\begin{pf}
The last statement follows easily from the other ones
and Remark \ref{rem-1}.
By definition, if $q_i$ is defective, then
$D\cdot E_i\ge -(B\cdot E_i)/2>0$.
Conversely, $D\cdot E_i\le0$ is equivalent to $(B+2K_Y)\cdot E_i\ge0$,
that holds, if $q_i$ is not defective, by Lemma 4.3 in \cite{ca}.
This proves the first statement.

Note that if there is a point $q_j>^1q_i$ with
$\tilde\alpha_j=\tilde\alpha_i$, then $q_j$ is the unique
proximate point to $q_i$.

Suppose that $q_i>^1q_k$ and $q_k$ is the only point which $q_i$ is
proximate to.
Then $\alpha_i=\tilde\alpha_i+\veps_k-\veps_i$ by \fref{alfa'}.
Let $q_j$ be an infinitely near point of order one to $q_i$.
If $q_j$ is proximate only to $q_i$,
then $\alpha_j=\tilde\alpha_j+\veps_i-\veps_j$,
thus $\alpha_j\ge\alpha_i+2$ if and only if
\[
(\tilde\alpha_i-\tilde\alpha_j)+\veps_k+\veps_j+2\le2\veps_i,
\]
which (recalling that $\tilde\alpha_i\ge\tilde\alpha_j$ because
$\tilde B\cdot E_i\ge0$) holds only if
$\tilde\alpha_i-\tilde\alpha_j=\veps_k=\veps_j=0$ and $\veps_i=1$,
that is case (i).

If $q_j$ is proximate also to $q_k$,
then $\alpha_i=\tilde\alpha_j+\veps_i+\veps_k-\veps_j$,
hence $\alpha_j\ge\alpha_i+2$ if and only if
$(\tilde\alpha_i-\tilde\alpha_j)+\veps_j+2\le2\veps_i$,
that is either case (i) or (ii) depending on the value of $\veps_k$.

This concludes the proof in case $q_i$ is proximate to only one point.
One may proceed similarly for the other configurations of $q_j>^1q_i$,
namely if $q_i$ is not infinitely near to any point
or if $q_i$ is proximate to more than one point.
\qed
\end{pf}

Both of the cases (i) and (ii) of Lemma \ref{Ldefective} may occur, as the
point $q_1$ in Example \ref{ex1} and the point $q_3$ in
Example \ref{ex3} respectively show.

We remark that if $q_i$ is defective and $q_j$ is as above,
namely $q_j >^1q_i$ and $\alpha_j>\alpha_i$, then $q_j$
cannot be defective. However there may exist a defective
point $q_l$ with $q_l>^1q_j$ and $\alpha_l=\alpha_i$,
as $q_3, q_5, \ldots, q_{2k-1}$ in Example \ref{ex2}.

We say that a point $q_i$ is \emph{$1$-defective}, and we
write $\ndef(q_i)=1$, if $q_i$ is defective and there is no
defective point $q_j>q_i$ with $\alpha_j=\alpha_i$. Inductively,
we say that $q_i$ is \emph{$k$-defective}, and we write
$\ndef(q_i)=k$, if there exists a $(k-1)$-defective point
$q_j>^2q_i$ with $\alpha_j=\alpha_i$.

In Example \ref{ex2}, the point $q_1$ is $k$-defective.
We set $\ndef(q_i)=0$ if $q_i$ is not defective
and $\Def=\{i\,|\,\ndef(q_i)>0\}$.\newcit{defDef}
Thus $i\in\Def$ if and only if $q_i$ is defective.

Now we are ready to show what exactly happens to an element in
$|\CI_\Gamma (h)|$ at a defective point.
To simplify the notation, by ordering conveniently the blowing-ups,
we may and will assume that
if $q_j,q_k$ are defective, with $\alpha_j=\alpha_k$ and $q_k>^2q_j$,
then $k=j+2$ and $q_{j+2}>^1q_{j+1}>^1q_j$.

\begin{theorem}\label{TbarE}
The fixed part of $|K_Y+B/2|$ is exactly:
\begin{equation}\label{barE}
\bar E=\sum_{j\in\Def}\sum_{r=0}^{\ndef(q_j)-1} E_{j+r}.
\end{equation}
\end{theorem}

\begin{pf}
By Lemma \ref{D.Ei}, the non-defective points do not mind,
so we may focus only on what happens at a defective point.
Any $k$-defective point looks like the point $q_1$ of Example \ref{ex2},
thus we will assume that $q_1$ is $k$-defective and we will follow
the notation of that example.
Recall that, for every $i=1$, $\ldots$, $k$,
the point $q_{2i-1}$ is $(k-i+1)$-defective,
$\alpha_{2i}/2-1=g$ and $\alpha_{2i-1}/2-1=g-1$.

We claim that the fixed part of $|K_Y+B/2|$ is exactly:
\begin{equation}\label{bbarE}
\bar E=\sum_{l=0}^{k-1}\sum_{i=1}^{k-l}E_{2i+l-1}
=\sum_{j\in\Def}\sum_{r=0}^{\ndef(q_j)-1} E_{j+r},
\end{equation}
and formula \fref{barE} clearly follows.
Note that \fref{bbarE} means that the general element
of $|{\cal I}_\Gamma(h)|$ has multiplicity $g$ at $q_1$, $\ldots$, $q_k$
and $g-1$ at $q_{k+1}$, $\ldots$, $q_{2k}$,
giving $kg^2$ adjunction conditions as expected.

Now we prove our claim.
By Lemma \ref{D.Ei}, the exceptional curve $E_{2i-1}$ is a fixed
component of $|K_Y+B/2|$, for $i=1$, $\ldots$, $k$. Then
\[
D+\sum_{i=1}^{k}E_{2i-1}
=\sum_{j=1}^{k}\left(gE_{2j-1}^*+(g-1)E_{2j}^*\right)
\]
meets positively $E_{2j}$ for $j=1,\ldots,k-1$, which therefore
are fixed components of $|K_Y+B/2|$ too. Now
\[
D+\sum_{i=1}^{k}E_{2i-1}+\sum_{j=1}^{k-1}E_{2j}
=gE_1^*+
\sum_{l=1}^{k-1}\left(gE_{2l}^*+(g-1)E_{2l+1}^*\right)
+(g-1)E_{2k}^*
\]
meets again positively $E_{2i-1}$, for $i=2$, $\ldots$, $k$
(but not $E_1$ and $E_{2k}$).

Going on in this way, by induction on
$k$, it follows that $|K_Y+B/2|$ contains \fref{bbarE}.
On the other hand,
$D+\bar E=\sum_{i=1}^k(gE_i^*+(g-1)E_{k+i}^*)$ does
not meet positively anyone of the $E_i$'s, thus the fixed part
of $|K_Y+B/2|$ is exactly $\bar E$, by Lemma \ref{Lultimo},
and our claim is proved.
\qed
\end{pf}

We remark that the previous theorem gives an alternative proof
of Corollary \ref{adjuncond}, independent from De Franchis' Theorem.

Finally we want to compute the number of conditions that the
singularity $p\in X_0$ imposes to {\em pluricanonical} systems.
The plurigenera of $X$ are:
\[
P_m(X)=h^0(X,mK_X)=h^0(Y,mK_Y+mB/2)+h^0(Y,mK_Y+(m-1)B/2).
\]
Riemann--Roch Theorem and \fref{as} imply that:
\[
h^0(mK_{X_0})-h^0(mK_X)+h^1(mK_X)
=\sum_{i=1}^n\frac{2(m^2-m)(\alpha_i-2)^2+\alpha_i^2-2\alpha_i}{8}
.
\]

\begin{theorem}\label{lastthm}
The number of conditions $h^0(mK_{X_0})-h^0(mK_X)$ that the
singularity $p\in X_0$ imposes to the $m$-canonical system are:
\begin{equation}\label{acmm}
\sum_{i=1}^n\frac{2(m^2-m)(\alpha_i-2)^2+\alpha_i^2-2\alpha_i}{8}-\frac{dm(m-1)}{2}
\end{equation}
where $d:=\sharp\Def$ is the number of defective points.
\end{theorem}

\begin{pf}
We shall show that $h^1(mK_X)=dm(m-1)/2$. Since we are dealing
with local questions, we may assume that the $(-1)$-curves of $X$
are contained in $\tau^{-1}(p)$. Recall that these $(-1)$-curves
are disjoint and there are exactly $d$ of them. Let
$\tau':X\to\bar X$ be their contraction (see section
\ref{minimalres}), then $\bar X$ is a minimal surface of general
type and we may assume that $h^0(K_{\bar X}) \gg 0$
(because $h$ is arbitrarily high). With no loss of
generality, we may also assume that $\tau'$ is the blowing-down of
just a $(-1)$-curve $F_i$. Thus it suffices to show that, under
these assumptions, $h^1(mK_X)=m(m-1)/2$. By Serre duality,
$h^1(mK_X)=h^1(-(m-1)K_X)$. Let $C$ be a curve in $|(m-1)K_X|$.
Clearly $C=\tau^{\prime*}(C_0)+(m-1)F_i$, where
$C_0\in|(m-1)K_{\bar X}|$.
It is well-known that $C_0$ and $\tau^{\prime*}(C_0)$ are 1-connected
(see \cite[Proposition 6.1]{bpv} and \cite[\S1]{fr2}),
therefore $h^0(\CO_{\tau^{\prime*}(C_0)})=1$.
Since $q(X)=h^1(\CO_X)=0$, the
exact sequence of sheaves
$
0\to \CO_X(-C) \to \CO_X \to \CO_C \to 0
$
implies that $h^1(\CO_X(-C))=h^0(\CO_{(m-1)F_i})$. Finally one
easily checks that $h^0(\CO_{(m-1)F_i})=m(m-1)/2$. \qed
\end{pf}

We remark that, if $h$ is not assumed to be arbitrarily high,
$\bar X$ may not be of general type and one should compute, or estimate,
$h^1(mK_X)$.

As we did for $|K_Y+B/2|$,
we want to determine the fixed components of $|mK_Y+mB/2|$ and
$|mK_Y+(m-1)B/2|$.
After having ordered the blowing-ups as explained just before Theorem \ref{TbarE},
we are ready to prove the following:

\begin{theorem}\label{verylastthm}
The fixed part of $|mK_Y+\bar mB/2|$, for $\bar m=m$
or $\bar m=m-1$, is exactly:
\begin{equation}
\left[\frac{\bar m}{2}\right]\sum_{j\in\Def}E_j+(\bar m\bmod 2)\bar E,
\end{equation}
where $\bar E$ is \fref{barE}, $[\bar m/2]$ is the largest integer
smaller than or equal to $\bar m/2$ and $\bar m \bmod 2=\bar m-2[\bar m/2]\in\{ 0, 1 \}$.
\qed
\end{theorem}

\begin{pf}
As in the proof of Theorem \ref{TbarE},
it suffices to understand what happens at a defective point,
thus we assume that $q_1$ is a $k$-defective point as in Example \ref{ex2}.
Let us set
\[
\tilde E=\sum_{i=1}^k E_{2i-1}=\sum_{j\in\Def} E_j.
\]
Suppose that $\bar m=m$.
If $m$ is even, then
\[
mD=\sum_{i=1}^{k}\left(m(g-1)E_{2i-1}^*+mgE_{2i}^*\right),
\]
hence $mD$ meets positively $E_{2i-1}$ for every $i=1,\ldots, k$.
Moreover $(mD+j\tilde E)\cdot E_{2i-1}>0$
for every $i=1,\ldots, k$ and $j=1,\ldots, m/2-1$.
This means that $m\tilde E/2$ is a fixed component of $|mK_Y+mB/2|$.
Then
\[
mD+\frac{m}{2}\tilde E=\sum_{i=1}^{2k}\left(mg-\frac{m}{2}\right)E_{i}
\]
does not meet positively anyone of the $E_i$'s, therefore
the fixed part of $|mK_Y+mB/2|$ is exactly $m\tilde E/2$,
by Lemma \ref{Lultimo}.

If $m$ is odd, following the same argument, one easily shows that the fixed
part of $|mK_Y+mB/2|$ is $(m-1)\tilde E/2+\bar E$.

One proceeds similarly in the case that $\bar m=m-1$.
Indeed, one can prove that the fixed part of $|mK_Y+(m-1)B/2|$
is $(m-1)\tilde E/2$ (resp.\ $m\tilde E/2+\bar E$) if $m$ is odd
(resp.\ if $m$ is even).
\qed
\end{pf}

Note that Theorem \ref{lastthm} can be proved also as
corollary of Theorem \ref{verylastthm}.

\end{document}